\documentclass{article}

\PassOptionsToPackage{numbers, square}{natbib}

\usepackage[final]{neurips_2022}

\title{Gradient Descent Is Optimal Under Lower Restricted Secant Inequality And Upper Error Bound}
\usepackage{times}
\author{%
  Charles Guille-Escuret \\
  Mila, Université de Montréal \\
  \texttt{guillech@mila.quebec} \\
  \And
  Adam Ibrahim \\
  Mila, Université de Montréal \\
  \texttt{adam.ibrahim@mila.quebec} \\
  \AND
  Baptiste Goujaud \\
  Ecole Polytechnique \\
  \texttt{baptiste.goujaud@gmail.com} \\
  \And
  Ioannis Mitliagkas \\
  Mila, Université de Montréal,\\ Canada CIFAR AI chair \\
  \texttt{ioannis@mila.quebec} \\
}
\usepackage[utf8]{inputenc} 
\usepackage[T1]{fontenc}    
\usepackage{hyperref}       
\usepackage{url}            
\usepackage{booktabs}       
\usepackage{amsfonts}       
\usepackage{nicefrac}       
\usepackage{microtype}      
\usepackage{faktor}
\usepackage{hhline}

\usepackage{amsmath}
\usepackage{amssymb}
\usepackage{wrapfig}

\usepackage{tikz}
\usepackage{mdframed}
\usepackage{subfigure}
\usepackage{float}

\mdfdefinestyle{myenvs}{
  nobreak=true, 
}

\newmdtheoremenv[style=myenvs]{theo}{Theorem}
\newmdtheoremenv[style=myenvs]{coro}{Corollary}
\newmdtheoremenv[style=myenvs]{lem}{Lemma}
\newmdtheoremenv[style=myenvs]{pro}{Prop}

\newtheorem{Th}{Theorem}[section]

\newtheorem{Def}[Th]{Definition}

\newtheorem{Rem}[Th]{Remark}
\newtheorem{?}[Th]{Problem}

\def\F{\mathcal{F}}

\def\EB{\operatorname{EB}}
\def\RSI{\operatorname{RSI}}

\def\R{\mathbb{R}}

\newcommand{\norm}[1]{\left\Vert#1\right\Vert_2}
\newcommand{\lrpar}[1]{\left( #1 \right)}

\newcommand{\braket}[2]{\left\langle #1 \ \middle| \ #2 \right\rangle}

\def\le{\left}
\def\ri{\right}
\def\eps{\epsilon}

\begin{document}

\maketitle

\begin{abstract}%
    The study of first-order optimization is sensitive to the assumptions made on the objective functions. 
These assumptions induce complexity classes which play a key role in worst-case analysis, including the fundamental concept of algorithm optimality.
Recent work argues that strong convexity and smoothness---popular assumptions in literature---lead to a pathological definition of the condition number \citep{guille2021study}.
Motivated by this result, we focus on the class of functions satisfying a lower restricted secant inequality and an upper error bound.
On top of being robust to the aforementioned pathological behavior and including some non-convex functions, this pair of conditions displays interesting geometrical properties. In particular, the necessary and sufficient conditions to interpolate a set of points and their gradients within the class can be separated into simple conditions on each sampled gradient. This allows the performance estimation problem (PEP, \cite{drori2012performance, taylor2016smooth}) to be solved analytically, leading to a lower bound on the convergence rate that proves gradient descent to be exactly optimal on this class of functions among all first-order algorithms.

\end{abstract}


\section{Introduction}

    \label{introduction}

The typical framework to study convergence properties of first-order algorithms in the context of machine learning is to first establish a class of objective functions to optimize through assumptions usually bound to a constant, such as $L$-smoothness and $\mu$-strong convexity. A tuning prescription of an algorithm is then made based on the constants, (e.g. step size $\alpha=\frac{2}{\mu+L}$ in the case of the gradient descent method on smooth and strongly convex functions), and finally a worst-case convergence rate can be derived for this algorithm when using this tuning prescription. In some cases, a lower bound on the achievable worst-case convergence rate can also be derived, leading to the theoretical optimality of an algorithm on the considered class of function, for instance the Nesterov accelerated gradient method \citep{Nesterov} is known to be optimal up to a constant on strongly convex and smooth functions.

However in \citep{guille2021study}, the authors establish that such framework and its derived results are very sensitive to the choice of assumptions, and that strong convexity and smoothness can exhibit pathological behaviors leading to conservative tuning and arbitrarily sub-optimal convergence rates, even when the resulting algorithm achieves theoretical worst-case optimality on this class of functions. Furthermore, they propose a set of more robust alternative conditions. In this work, we focus on a specific pair of such alternative conditions : lower restricted secant inequality ($\RSI^-$) and upper error bounds ($\EB^+$). Our main contribution is to show that the gradient descent (GD) method with a certain tuning is exactly optimal on the classes of objective functions induced by these conditions, confirming that optimality results are highly sensitive to the choice of conditions. Another consequence is that no algorithm can accelerate on this class of functions, implying that additional assumptions are required to explain the practical efficiency of accelerated methods.

\paragraph{Notation}
Let $\mathcal{F}$ the set of differentiable functions from $\R^d$ to $\R$ that admit a convex set of global minima $X_f^*$. We focus on the problem of optimizing a function $f\in \F$, i.e. finding $x\in X_f^*$. For any $x\in\R^d$, we denote $x_f^*$ the orthogonal projection of $x$ on $X_f^*$. By abuse of notation, when the context is not ambiguous, we will simply denote $x_f^*$ as $x^*$.

We call gradient descent (GD) the standard optimization algorithm based on the following update, 
where $\alpha$ is the step size : 
$$x_{i+1} = x_i - \alpha \nabla f(x_i)$$ 

We call \emph{first-order algorithm} all $\mathcal{A}$ that consider past iterates, function values and gradients and output a next iterate. Formally, $\mathcal{A}$ can be seen as a sequence of functions $\le\{\mathcal{A}_n\mid n\in\mathbb{N}\ri\}$ such that for any $n\in\mathbb{N}$, $\mathcal{A}_n$ is a function defined on $\le(\R^d\times \R \times \R^d\ri)^{n+1}$ with values in $\R^d$. Under that formalism, applying $\mathcal{A}$ to optimize an objective function $f$ starting in $x_0\in\R^d$ generates a sequence of iterates $(x_i)_i$ such that $\forall i, x_{i+1} = \mathcal{A}_{i}\le(\le(x_j,f(x_j),\nabla f(x_j)\ri)_{j\leq i}\ri)$. Note that we do not require the iterates of $\mathcal{A}$ to lie within the span of observed gradients as it is often the case in the literature.

\paragraph{Outline} Section \ref{sec:RSIEB} introduces $\RSI^-$ and $\EB^+$ and provides some basic properties and motivation. Section \ref{related_work} discusses related works in the literature. Section \ref{interpolation} defines and establishes the necessary and sufficient interpolation conditions for $\RSI^-$ and $\EB^+$, which is a key element of our analysis. Section \ref{lower_bound} proves the lower bound on the convergence rate of first-order algorithms under $\RSI^-$ and $\EB^+$, and finally Section \ref{conclusion} concludes our work. Detailed proofs are provided in the appendix.

\section{Lower restricted secant inequality and upper error bounds}

    \label{sec:RSIEB}

We now define $\RSI^-$ and $\EB^+$ and discuss some basic properties. 
\begin{Def}[Lower restricted secant inequality]
Let $f\in\F$ and $\mu>0$
$$f\in \RSI^-(\mu) \Leftrightarrow \forall x\in\R^d, \le\langle \nabla f(x) \mid x-x^*\ri\rangle \geq \mu\le\|x-x^*\ri\|_2^2$$
\end{Def}

Intuitively, $\RSI^-(\mu)$ enforces that the further $x$ is from $X_f^*$, the stronger the gradient of $f$ in $x$ must be in the opposite direction of $X_f^*$.

\begin{Rem}
$\RSI^-(\mu)$ includes non-convex functions. However, it prevents flat landscapes outside of $X_f^*$, and requires $f$ to increase at least quadratically with the distance to $X_f^*$, as established in \citep{guille2021study} : $$f\in \RSI^-(\mu)\Rightarrow\forall x\in\R^d, f(x)-f^*\geq \frac{\mu}{2}\le\|x-x^*\ri\|_2^2$$
\end{Rem}

\begin{Def}[Upper Error Bounds]
Let $f\in\F$ and $L>0$
$$f\in \EB^+(L) \Leftrightarrow \forall x\in\R^d, \le\| \nabla f(x)\ri\|_2 \leq L\le\|x-x^*\ri\|_2$$
\end{Def}

$\EB^+(L)$ thus enforces that the gradient of $f$ is controlled by the distance to $X_f^*$. 
\begin{Rem}
\label{rem_intro}
$L$-smoothness implies $\EB^+(L)$ and $\mu$-strong convexity implies $\RSI^-(\mu)$. However, one must be careful before claiming that $\RSI^-$ and $\EB^+$ are respectively weaker than strong convexity and smoothness : for $\mu_0>\mu_1$, $\RSI^-(\mu_0)$ is neither stronger or weaker than $\mu_1$-strong convexity, and for $L_0<L_1$, $\EB^+(L_0)$ is neither stronger or weaker than $L_1$-smoothness.
\end{Rem}
Therefore, even when the objective function is smooth and strongly convex, considering convergence results under $\RSI^-$ and $\EB^+$ is relevant, as we might obtain better constants $\mu$ and $L$. Many convergence results depend of the condition number $\kappa=\frac{L}{\mu}$. Better constants leads to a better condition number, and thus a potentially better convergence rate (including when the dependence in the condition number $\kappa$ is the same or worse). As a consequence, machine learning problems with strongly convex and smooth objective functions are all potential applications of results under $\RSI^-$ and $\EB^+$, provided we can obtain better constants under these conditions.

\begin{pro}[convergence rate of GD under $\RSI^-$ and $\EB^+$]
\label{prop_intro}
Let $f\in \RSI^-(\mu)\cap \EB^+(L)$. Then gradient descent with learning rate $\alpha=\frac{\mu}{L^2}$ on $f$ will guarantee the following convergence rate: 
\begin{equation}
\label{prop:eq1}
\le\|x_{i}-x_{i}^*\ri\|_2^2\leq \le(1-\frac{\mu^2}{L^2}\ri)^i\le\|x_0-x_0^*\ri\|_2^2
\end{equation}
Moreover, when the learning rate is set to $\alpha_0=\frac{1}{2\mu}$ on the first step, and $\alpha=\frac{\mu}{L^2}$ on every other step, gradient descent guarantees the following convergence rate:
\begin{equation}
\label{prop:eq2}
\le\|x_{i}-x_{i}^*\ri\|_2^2\leq \frac{\le\|\nabla f(x_0)\ri\|_2^2}{4\mu^2}\le(1-\frac{\mu^2}{L^2}\ri)^{i-1}
\end{equation}
\end{pro}
\textit{Proof.} When $\alpha=\frac{\mu}{L^2}$, we have
\begin{equation}
\begin{aligned}
    \le\|x_{i+1}-x_{i+1}^*\ri\|_2^2 & \leq \le\|x_{i+1}-x_i^*\ri\|_2^2 \\
        & = \le\|x_i-\alpha \nabla f(x_i)-x_i^*\ri\|_2^2 \\
        & = \le\|x_i-x_i^*\ri\|_2^2 - 2 \alpha \le\langle \nabla f(x_i)\mid x_i-x_i^*\ri\rangle + \alpha^2\le\|\nabla f(x_i)\ri\|_2^2 \\
        & \leq\le(1-2\alpha\mu+L^2\alpha^2\ri)\le\|x_i-x_i^*\ri\|_2^2 \\
        & = \le(1-\frac{\mu^2}{L^2}\ri)\le\|x_i-x_i^*\ri\|_2^2 \,,
\end{aligned}
\end{equation}
which proves \eqref{prop:eq1}. To prove \eqref{prop:eq2}, we simply note that when using $\alpha=\frac{1}{2\mu}$, we have:
\begin{equation}
\begin{aligned}
    \le\|x_{1}-x_{1}^*\ri\|_2^2 & \leq \le\|x_0-x_0^*\ri\|_2^2 - 2 \alpha \le\langle \nabla f(x_0)\mid x_0-x_0^*\ri\rangle + \alpha^2\le\|\nabla f(x_0)\ri\|_2^2 \\
        & \leq \frac{1}{\mu}\le\langle \nabla f(x_0)\mid x_0-x_0^*\ri\rangle-\frac{1}{\mu}\le\langle \nabla f(x_0)\mid x_0-x_0^*\ri\rangle +\frac{\le\|\nabla f(x_0)\ri\|_2^2}{4\mu^2} \\
        & = \frac{\le\|\nabla f(x_0)\ri\|_2^2}{4\mu^2} \,.
\end{aligned}
\end{equation}

$\hfill\blacksquare$

Interestingly, $\RSI^-$ and $\EB^+$ are direct bounds on the two additional terms obtained by developing $\le\|x_i-\alpha\nabla f(x_i)-x_i^*\ri\|_2^2$, leading to an extremely simple proof. On the intuitive level, $\RSI^-$ lower bounds the gain from stepping in the direction of $X_f^*$, while $\EB^+$ upper bounds the error coming from the component of the gradient orthogonal to that direction.

\begin{Rem}
The literature gives a worst case convergence rate of gradient descent on $\mu$-strongly convex and $L$-smooth functions of $\le\|x_{i}-x_{i}^*\ri\|_2^2\leq \le(1-\frac{2\mu}{\mu+L}\ri)^{2i}\le\|x_0-x_0^*\ri\|_2^2$,
using the step size $\alpha=\frac{2}{\mu+L}$ \citep{nesterov2003introductory,polyak1987introduction}. While this rate is better for a fixed $\mu$ and $L$, we again emphasize that the constants $\mu$ and $L$ may be very different depending on the chosen conditions, and thus these two rates can not directly be compared. 
\end{Rem}

\textbf{Motivation:} we empirically verify in Appendix \ref{appendix:rsiebexp} that the optimization path of a ResNet18~\cite{https://doi.org/10.48550/arxiv.1512.03385} trained for classification on CIFAR10~\cite{Krizhevsky09learningmultiple} verifies the interpolation conditions of $RSI^-\cap EB^+$, which are introduced in Section \ref{interpolation}. This result guarantees the existence of a function in $RSI^-\cap EB^+$ which exactly interpolates the observed gradients, and thus all convergence guarantees of $RSI^-\cap EB^+$ naturally apply in this practical deep learning setting. $RSI^-\cap EB^+$ therefore provides linear convergence guarantees while being empirically applicable to neural networks (although with respects to a local minima), an impressive feat given the highly non-convex nature of neural networks loss functions.

\section{Related Work}

    \label{related_work}
Throughout the literature, many choices of assumptions have been used to study first-order optimization. Most assumptions fall into one of two categories : \textit{lower conditions} and \textit{upper conditions} that respectively take the form of a lower and an upper bound on properties of the objective function. For instance, strong convexity lower bounds the curvature of the objective function and is thus a \textit{lower condition}. Similarly, smoothness is an \textit{upper condition}. 

Lower conditions have been the most extensively studied assumptions, such as Polyak-{\L}ojasiewicz \citep{polyak_gradient_1963}, \emph{local-quasi-convexity} \citep{NIPS2015_5718}, \emph{weak quasi-convexity} \citep{hardt2016gradient}, \emph{quadratic growth} \citep{anitescu2000degenerate,bonnans1993second, ioffe1994sensitivity}, \emph{Kurdyka-{\L}ojasiewicz} \citep{kurdyka1998gradients,bolte2008characterizations}, \emph{optimal strong convexity} \citep{liu2014asynchronous,ma2015linear,gong2014linear}, \emph{weak strong convexity} \cite{DBLP:journals/corr/KarimiNS16,necoara2016linear}, \emph{error bounds} \citep{luo_error_1993}. Some recent works have explored the relations between these lower conditions \citep{DBLP:journals/corr/KarimiNS16,zhang2017restricted}. In this work we focus on the \emph{restricted secant inequality} ($\RSI^-$) (which we denote as \emph{lower restricted secant inequality} to differentiate it from its upper bound equivalent) which was introduced in \citep{zhang2013gradient}, and has been used (along with its convex extension \emph{restricted strong convexity}) in many recent theoretical derivations of linear convergence rates \citep{yi2019exponential,Schpfer2016LinearCO,Yuan2016OnTC}.

On the contrary, because most machine learning objective functions are naturally smooth, fewer works have explored alternatives to smoothness. However as discussed in Remark \ref{rem_intro}, it is still relevant to study these alternatives on smooth objectives due to potentially better conditioning. The most notable ones in the literature are \emph{local smoothness} \citep{NIPS2015_5718}, \emph{restricted smoothness} \citep{agarwal2011fast}, \emph{relative smoothness} \citep{lu2016relativelysmooth, hanzely2018accelerated, zhou2019simple}, \emph{weak-smoothness} \citep{hardt2016gradient}, \emph{expected smoothness} \cite{gower2020sgd}. In \citep{guille2021study}, the authors argue that \textit{lower conditions} can be naturally translated into equivalent \textit{upper conditions} by changing the lower bound into an upper bound, and vice-versa. Subsequently, they introduce a set of \textit{upper} equivalent to existing \textit{lower conditions}, such as \emph{upper error bounds} $\EB^+$, the natural \textit{upper} equivalent to \emph{error bounds} from \citep{luo_error_1993}. Throughout this work, we focus on $\EB^+$ as an \textit{upper condition}.

Finally, a key to our analysis are the necessary and sufficient interpolation conditions of $\RSI^-(\mu)\cap \EB^+(L)$ (Section \ref{interpolation}). The search of these conditions has been largely motivated by Performance Estimation Problems (PEP), introduced in \citep{taylor2016smooth} and with many recent successful applications \citep{hu2017dissipativity, Taylor2017ConvexIA,2019, taylor2020exact,2017}. PEP is a framework for computer-assisted worst-case convergence analysis, that only requires necessary and sufficient interpolation conditions for the considered class of objective functions. In the case of $\RSI^-\cap \EB^+$ however, we found the interpolation conditions to be \emph{independent} (see Section \ref{interpolation}), which made the worst-case analysis directly solvable analytically.

\section{Interpolation conditions}

    \label{interpolation}

In this section we provide and discuss the necessary and sufficient interpolation conditions for $\RSI^-\cap \EB^+$. Their importance stems from the framework used in PEP and in this work to analyze worst-case convergence. Given a first-order optimization algorithm $\mathcal{A}$ we want to find the slowest convergence rate of $\mathcal{A}$ among all $f$ within a class of objective functions $\mathcal{C}$ and starting point $x_0\in\R^d\setminus X_f^*$. That is equivalent to solving the following optimization problem at any step number $n$ : 
\begin{equation}
\label{pep1}
\begin{aligned}
\min_{f\in\mathcal{C},(x_i)_{i\leq n}\in\left(\R^d\right)^{n+1}} \quad & \frac{\le\|x_0-x_0^*\ri\|_2}{\le\|x_n-x_n^*\ri\|_2}\\
\textrm{s.t.} \quad & \forall i\leq n-1, x_{i+1} = \mathcal{A}\le(\le(x_0,f(x_0),\nabla f(x_0)\ri), ..., \le(x_i, f(x_i),\nabla f(x_i)\ri)\ri)\\
\end{aligned}
\end{equation}
Directly searching for $f$ in the functional space is generally intractable, however if we can explicitly find the set $\mathcal{G}$ of all families $(x_i, f_i, g_i)_i$ such that $\exists f\in\mathcal{C}, \forall i, \nabla f(x_i)=g_i$ and $f(x_i)=f_i$  (we say that $f$ interpolates $(x_i,f_i,g_i)_i$), then problem \eqref{pep1} can be reduced to:
\begin{equation}
\label{pep2}
\begin{aligned}
\min_{(x_i, g_i,f_i)_{i\leq n}\in(\R^d \times \R^d \times \R)^{n+1}} \quad & \frac{\le\|x_0-x_0^*\ri\|_2}{\le\|x_n-x_n^*\ri\|_2}\\
\textrm{s.t.} \quad & \forall i\leq n-1, x_{i+1} = \mathcal{A}\le(\le(x_0,f_0,g_0\ri), ..., \le(x_i, f_i, g_i\ri)\ri)\\
\textrm{and} \quad & (x_i,f_i,g_i)_{i\leq n} \in\mathcal{G} \\
\end{aligned}
\end{equation}
In many cases (see Section \ref{related_work}), problem \eqref{pep2} is tractable and becomes a very powerful analysis tool providing lower bounds, upper bounds, and optimal tuning for different types of algorithms and assumptions used. A crucial and difficult component of this analysis is to formulate the \textit{interpolation conditions}, that is the necessary and sufficient conditions for a family $(x_i,f_i,g_i)_i$ to belong in $\mathcal{G}$. Driven by these motivations, we now establish the interpolation conditions for $\RSI^-\cap \EB^+$. 

In Theorem \ref{interpolation_th}, we introduce the necessary and sufficient conditions to interpolate a family $(x_i,g_i)$, without considering the function values $(f_i)_i$. This theorem could be used to find the worst case convergence rate over all first-order algorithms that ignore function values. However, in Corollary \ref{interp_cor}, we deduce from Theorem \ref{interpolation_th} sufficient (but not necessary) conditions to interpolate a family $(x_i,f_i,g_i)$. These conditions allow us in Section \ref{lower_bound} to find a lower bound on the worst-case convergence rate for \emph{all} first-order algorithms (including algorithms that have access to function values information), which we know is tight thanks to Prop \ref{prop_intro}.

\begin{theo}[Interpolation conditions]
\label{interpolation_th}
Let $(x_i, g_i)_{i\leq n} \in \le(\R^d\times\R^d\ri)^{n+1}$, such that the $x_i$ are separate points.

Then, $\forall \mu, L>0$:
$$\exists f\in \RSI^-(\mu)\cap \EB^+(L),\: s.t.\: \forall i, \: \nabla f(x_i)=g_i$$
$$\Updownarrow$$
$$\exists X^*\subseteq \mathbb{R}^d\ convex,\; s.t.\; \forall i,$$
\begin{equation}
    \label{interp}
    \begin{aligned}
        \le\|g_i\ri\|_2 \leq L\le\|x_i-x_i^*\ri\|_2\quad \text{and}\quad \left\langle g_i \mid x_i-x_i^* \right\rangle \geq \mu\le\|x_i-x_i^*\ri\|_2^2,
    \end{aligned}
\end{equation}

where $x_i^*$ is the orthogonal projection of $x_i$ onto $X^*$.
\end{theo}
\textit{Proof.} In order to preserve concision and clarity, we will only present the broad outline of the proof here. For a complete technical proof, see Appendix \ref{proof_interpolation}.

The direct implication is trivial as it is a direct application of $\RSI^-(\mu)$ and $\EB^+(L)$ definitions to $f$ in $(x_i)_i$. For the reverse implication, we will construct a function $f_{\eps,\beta}$ that interpolates $(x_i,g_i)_i$. The function $f_{\eps,\beta}$ is a quadratic everywhere except in the spheres of radius $\eps$ around each $x_i$, $\eps$ being small enough for these spheres to never intersect. Inside the sphere of radius $\eps$ around a given $x_i$, $f_{\eps,\beta}$ will be perturbed by adding a term $\lambda(\le\|x-x_i\ri\|_2)h(x)$ where $h$ is affine in $x$, and $\lambda$ is a scaling term so that $\lambda(\eps)=0$ at the border of the sphere, and $\lambda(0)=1$ in its center $x_i$.

The key is to find a function $\lambda$ that preserves the properties of $\RSI^-(\mu)$ and $\EB^+(L)$. We use 
\begin{equation}
    \lambda_{\eps,\beta}(u) = \frac{1+\cos\le(\pi\frac{u^\beta}{\eps^\beta}\ri)}{2}
\end{equation}
And our construction $f_{\eps,\beta}$ is given by:
\begin{gather}
f_{\eps,\beta}(x) = \begin{cases}
\frac{\mu+L}{4}\le\|x-x^*\ri\|_2^2 & \text{if}\; \forall i, \le\|x-x_i\ri\|_2 \geq \eps\\
\frac{\mu+L}{4}\le\|x-x^*\ri\|_2^2 + \lambda_{\eps,\beta}\le(\le\|x-x_i\ri\|_2\ri)\le\langle g_i-\frac{\mu+L}{2}(x_i-x_i^*)\mid x-x_i\ri\rangle & \text{if}\; \exists i, \le\|x-x_i\ri\|_2 < \eps
\end{cases}
\end{gather}
The rest of the proof is to use the Taylor expansions of $f_{\eps,\beta}$ to show that for sufficiently small $\eps$ and $\beta$, $f_{\eps,\beta}$ will belong in $\RSI^-(\mu)$ and $\EB^+(L)$ (see Appendix \ref{proof_interpolation}).

$\hfill\blacksquare$

\begin{coro}
\label{interp_cor}
Let $(x_i, f_i, g_i)_{i\leq n} \in \le(\R^d\times\R\times\R^d\ri)^{n+1}$, such that the $x_i$ are separate points. 

Then, $\forall \mu, L>0$:
$$\exists X^*\subseteq \mathbb{R}^d\ convex,\; s.t.\; \forall i,$$
\begin{equation}
    \begin{aligned}
     \le\|g_i\ri\|_2 & \leq L\le\|x_i-x_i^*\ri\|_2 \\
     \left\langle g_i \mid x_i-x_i^* \right\rangle & \geq \mu\le\|x_i-x_i^*\ri\|_2^2 \\
     f_i & = \frac{\mu+L}{4}\le\|x_i-x_i^*\ri\|_2^2\\
    \end{aligned}
\end{equation}

$$\Downarrow$$
$$\exists f\in \RSI^-(\mu)\cap \EB^+(L),\: s.t.\: \forall i, \: \nabla f(x_i)=g_i \quad and \quad f(x_i)=f_i,$$

where $x_i^*$ is the orthogonal projection of $x_i$ onto $X^*$.
\end{coro}
\textit{Proof.} We simply use the function $f_{\eps,\beta}$ from the proof of Theorem \ref{interpolation_th} and note that $\forall i, f_{\eps,\beta}(x_i)=\frac{\mu+L}{4}\le\|x_i-x_i^*\ri\|_2^2$.

$\hfill\blacksquare$

Theorem \ref{interpolation_th} and Corollary \ref{interp_cor} are the key elements to prove the optimality of gradient descent on $\RSI^-$ and $\EB^+$ among all first-order algorithms (see Section \ref{lower_bound}).

\begin{Rem}
\label{rem:indep}
The interpolation conditions in Theorem \ref{interpolation_th} are independent, in the sense that a family $(x_i,g_i)_i$ admits an interpolation in $\RSI^-(\mu)\cap \EB^+(L)$ if and only if each $\{(x_i,g_i)\}$ admits an interpolation in $\RSI^-(\mu)\cap \EB^+(L)$. Similarly, the sufficient interpolation conditions in Corollary \ref{interp_cor} are also independent.
\end{Rem}

This property of independent interpolation conditions drastically simplifies convergence analysis and is the main reason we are able to analytically derive a lower bound in Section \ref{lower_bound}. Indeed, given an interpolable family $(x_i,f_i,g_i)_{i\leq n}$ for a given set $X^*$, it is sufficient to show that $(x_{n+1},f_{n+1},g_{n+1})$ is interpolable with $X^*$ to prove that the entire family $(x_i,f_i,g_i)_{i\leq n+1}$ can be interpolated. It is thus simple to find the set of interpolable $(f_{n+1},g_{n+1})$ given $X^*$, $x_{n+1}$, and an interpolable $\le(x_i,f_i,g_i\ri)_{i\leq n}$.

\section{Lower bound on the convergence rate}

    \label{lower_bound}
In this Section we derive a lower bound on the convergence rate of first-order algorithms on $\RSI^-$ and $\EB^+$. This lower bounds applies under the assumption that the number of steps taken is smaller than the number of dimension $d$. This assumption is frequent in the literature (e.g. \cite{bubeck2015convex}) and not constraining for high-dimensional optimization. The observed gradients of the worst-case functions for optimal algorithms are typically orthogonal to one another (see \cite{2019}) which is not possible when the number of steps becomes larger than the dimension $d$. We conjecture that when not bounding the number of steps, it is possible to achieve an asymptotic rate in $O\le(2^{-\frac{n}{d}}\ri)$ which would be better than the usual rates obtained for very ill-conditioned functions, while having little to no practical uses due to the bad convergence properties on a lower number of steps.

We now introduce Lemma \ref{lemma_lower}, which is the cornerstone of the proof of Theorem \ref{lower_bound_th}:

\begin{lem}
\label{lemma_lower}
Let $\mu>0$ and $L>\mu$. Let $\alpha_0\in\le[\frac{\mu}{L^2},\max\le(\frac{\mu}{L^2},\frac{1}{2\mu}\ri)\ri]$. For any first-order optimization algorithm $\mathcal{A}$ and starting point $x_0\in\R^d$, there exists $(g_i)_{i\leq d-2}\in\mathbb{R}^{d}$, $(f_i)_{i\leq d-2}\in\R$ and  $\mathcal{S}_{d-2}\subseteq\mathcal{S}_{d-1}\subseteq\dots\subseteq\mathcal{S}_0\subseteq\R^{d}$ such that:
\begin{enumerate}
    \item $\forall i\leq d-2$, there exists a $(d-i-1)$-dimensional affine space $\mathcal{H}_i$ containing $\mathcal{S}_i$ and in which $\mathcal{S}_i$ is a $(d-i-2)$-sphere of radius $r_i=\sqrt{\frac{\alpha_0}{\mu}-\alpha_0^2}\le\|g_0\ri\|_2\le(1-\frac{\mu^2}{L^2}\ri)^{\frac{i}{2}}$ and center $c_i\in\mathcal{H}_i$.
    \item Let $(x_i)_i$ be the iterates generated by $\mathcal{A}$ starting from $x_0$ and reading gradients $(g_i)_i$ and function values $(f_i)_i$, then for any $i\leq d-2$ and any $x\in\mathcal{S}_i$, there exists a function $f$ in $\RSI^-(\mu)\cap \EB^+(L)$ minimized by $\{x\}$ that interpolates $(x_j,f_j,g_j)_{j\leq i}$.
\end{enumerate}
\end{lem}
\textit{Proof.} In order to preserve concision and clarity, we will only present the broad outline of the proof here. For a complete technical proof, see Appendix \ref{lowerbound_proof}.

We construct the sequence iteratively. For initialisation, we take any non-zero $g_0$, set $f_0=\frac{\mu+L}{4\mu}\alpha_0\le\|g_0\ri\|_2^2$, $c_0=x_0-\alpha_0 g_0$, and finally
$$\mathcal{S}_0 = \le\{x\in\R^d\middle| \le\langle x-c_0\mid g_0\ri\rangle = 0\ri\}\cap\le\{x\in \R^d\middle| \le\|x-c_0\ri\|_2=\sqrt{\frac{\alpha_0}{\mu}-\alpha_0^2}\le\|g_0\ri\|_2\ri\}$$
Then assuming we have a sequence $(f_j,g_j,\mathcal{S}_j)_{j\leq i<(d-2)}$ respecting the conditions of the Lemma, noting $\mathcal{H}_i$ the $d-i-1$ dimensional affine space in which $\mathcal{S}_i$ is a $(d-i-2)$ dimensional sphere, and $x_{i+1}$ the $(i+1)$-th iterate returned by $\mathcal{A}$. Let $h_{i+1}$ the orthogonal projection of $x_{i+1}$ into $\mathcal{H}_i$.

If $h_{i+1}\neq c_i$, let $v=\frac{(h_{i+1}-c_i)}{\le\|h_{i+1}-c_i\ri\|_2}$. If $h_{i+1}=c_i$, let $s\in\mathcal{S}_{i}$ and $v=\frac{(s-c_i)}{\le\|s-c_i\ri\|_2}$. 

We then construct:

$$c_{i+1}=c_i-\frac{\mu}{L}r_i v$$

$$f_{i+1}=\frac{\mu+L}{4}(\le\|x_{i+1}-c_{i+1}\ri\|_2^2+(1-\frac{\mu^2}{L^2})r_i^2)$$

$$g_{i+1}=L\frac{\le\|x_{i+1}-x^*\ri\|_2}{\le\|x_{i+1}-c_{i+1}\ri\|_2}\le(x_{i+1}-c_{i+1}\ri)$$

$$\mathcal{H}_{i+1}=\le\{x\in\mathcal{H}_i\mid \le\langle x-c_i\mid v\ri\rangle=-\frac{\mu}{L}r_i \ri\}$$

$$\mathcal{S}_{i+1}=\mathcal{S}_i\cap\mathcal{H}_{i+1}$$

We verify in Appendix \ref{lowerbound_proof} that this construction respects the properties of Lemma \ref{lemma_lower}.

$\hfill\blacksquare$

We can now introduce Theorem \ref{lower_bound_th} which gives us a lower bound on the worst-case convergence rate of any first-order algorithm on $\RSI^-(\mu)$ and $\EB^+(L)$.

\begin{theo}[Lower bound on $\RSI^-\cap \EB^+$]
\label{lower_bound_th}

Let $\mathcal{A}$ be any first-order algorithm on $\mathbb{R}^d$, $\mu>0$ and $L\geq\mu$. For any $x_0\in\mathbb{R}^d$, there exists $x^*\in\R^d$ and a function $f$ in $\RSI^-(\mu)\cap \EB^+(L)$ minimized by $X^*=\{x^*\}$ such that  
\begin{equation}
\label{eq:th1}
\forall i \leq d-1, \le\|x_i-x_i^*\ri\|_2^2 \geq \le(1-\frac{\mu^2}{L^2}\ri)^i\le\|x_0-x_0^*\ri\|_2^2
\end{equation}
Furthermore, if $\frac{L}{\mu}\geq \sqrt{2}$, there exists $h$ in $\RSI^-(\mu)\cap \EB^+(L)$ minimized by $X^*=\{x^*\}$ such that
\begin{equation}
\label{eq:th2}
\forall i \leq d-1, \le\|z_i-z_i^*\ri\|_2^2 \geq \frac{\le\|\nabla h(z_0)\ri\|_2^2}{4\mu^2}\le(1-\frac{\mu^2}{L^2}\ri)^{i-1}
\end{equation}
where $(x_i)$ (resp. $(z_i)$) is the trajectory obtained by applying $\mathcal{A}$ to $f$ (resp. $h$) starting in $x_0$.

\end{theo}

Note that since $X^*$ is a singleton, $\forall i, x_i^* = z_i^* = x^*$.

\textit{Proof.} If $L=\mu$, then the inequalities are trivial from the positivity of the norm. If $L>\mu$, let $(g_i,f_i, \mathcal{S}_i)_{i\leq (d-2)}$ be the sequence introduced in Lemma \ref{lemma_lower} for $\alpha_0\in\le[\frac{\mu}{L^2},\max\le(\frac{\mu}{L^2},\frac{1}{2\mu}\ri)\ri]$ . Let us note that for any $x\in\mathcal{S}_0$, $\le\|x_0-x\ri\|_2^2=\frac{\alpha_0\le\|g_0\ri\|_2^2}{\mu}$ (see initialisation in Appendix \ref{lowerbound_proof}). 

Let $i\in\le\{1,\dots,d-1\ri\}$. $\mathcal{S}_{i-1}$ has radius $r_{i-1}$, thus there exists $x^*\in\mathcal{S}_{i-1}$ such that $\le\|x_{i}-x^*\ri\|_2\geq r_{i-1}$, and thus :

\begin{equation}
\label{eq:th3}\le\|x_{i}-x^*\ri\|_2^2\geq r_{i-1}^2=\le(\frac{\alpha_0}{\mu}-\alpha_0^2\ri)\le\|g_0\ri\|_2^2\le(1-\frac{\mu^2}{L^2}\ri)^{i-1}
\end{equation}

When setting $\alpha_0=\frac{\mu}{L^2}$ and observing that $x^*\in\mathcal{S}_{i-1}\subseteq\mathcal{S}_0$ and thus $\le\|x_0-x^*\ri\|_2^2=\frac{\alpha_0\le\|g_0\ri\|_2^2}{\mu}$ in \eqref{eq:th3}, we obtain \eqref{eq:th1}. If $\frac{L}{\mu}\geq\sqrt{2}$, we set $\alpha_0=\frac{1}{2\mu}$ and \eqref{eq:th3} immediately yields \eqref{eq:th2}.

$\hfill\blacksquare$

\begin{Rem}
Since the lower bounds established in Theorem \ref{lower_bound_th} are exactly matched by the convergence guarantees of gradient descent (see Prop \ref{prop_intro}), these bounds are tight and gradient descent is exactly optimal on $\RSI^-(\mu)\cap\EB^+(L)$. This is a concrete example of the sensitivity of theoretical optimality to the choice of complexity classes.
\end{Rem}

\begin{subsection}{Discussion}
\label{sec:discussion}

The first bound presented in Theorem \ref{lower_bound_th} gives the optimal solution when trying to solve
\begin{equation}
\label{eq:rem1}
\min_{\mathcal{A}}\max_{f,x_0} \quad \frac{\le\|x_n-x_n^*\ri\|_2}{\le\|x_0-x_0^*\ri\|_2}
\end{equation} 
While the second bound gives the optimal solution when trying to solve 
\begin{equation}
\label{eq:rem2}
\min_{\mathcal{A}}\max_{f,x_0} \quad \frac{\le\|x_n-x_n^*\ri\|_2}{\le\|\nabla f(x_0)\ri\|_2}
\end{equation}  
For general smooth and convex functions, \eqref{eq:rem2} will not have a solution (for any $\mathcal{A}$, the quantity will not have a worst case upper bound), which is why \eqref{eq:rem1} has historically been the focus of optimization literature. However in practice, when \eqref{eq:rem2} admits a solution, as is the case for $\RSI^-(\mu)\cap\EB^+(L)$, it fits practical motivations better than \eqref{eq:rem1} : when starting from $x_0$ and observing an initial gradient $g_0$, the solution to \eqref{eq:rem2} is the one that will minimize $\le\|x_n-x_n^*\ri\|_2$ in the worst case. In comparison, for a fixed $\le\|\nabla f(x_0)\ri\|_2$, the solution to \eqref{eq:rem1} will be faster when $\le\|x_0-x_0^*\ri\|_2$ is small, and slower when $\le\|x_0-x_0^*\ri\|_2$ is large, leading to a slower worst-case convergence.

When $\frac{L}{\mu}>\sqrt{2}$, using a tuning of $\alpha_0=\frac{1}{2\mu}$ on the first step instead of $\frac{\mu}{L^2}$ leads to a worst-case convergence of $\le\|x_n-x_n^*\ri\|_2$ better by a constant $c=\frac{L^2}{2(L^2-\mu^2)}$. While such small constant factor is often considered not impactful, the number of steps required to make up for this constant factor is $n=\frac{-\log(2)}{\log\le(1-\frac{\mu^2}{L^2}\ri)}-1$, which yields $n\approx 68$ for $\frac{L}{\mu}=10$ and $n\approx 6930$ for $\frac{L}{\mu}=100$, and can thus become substantial on ill-conditioned functions.

\begin{figure}
\centering     
\subfigure[$\frac{L}{\mu}>\sqrt{2}$]{\label{fig:1}\includegraphics[width=39mm]{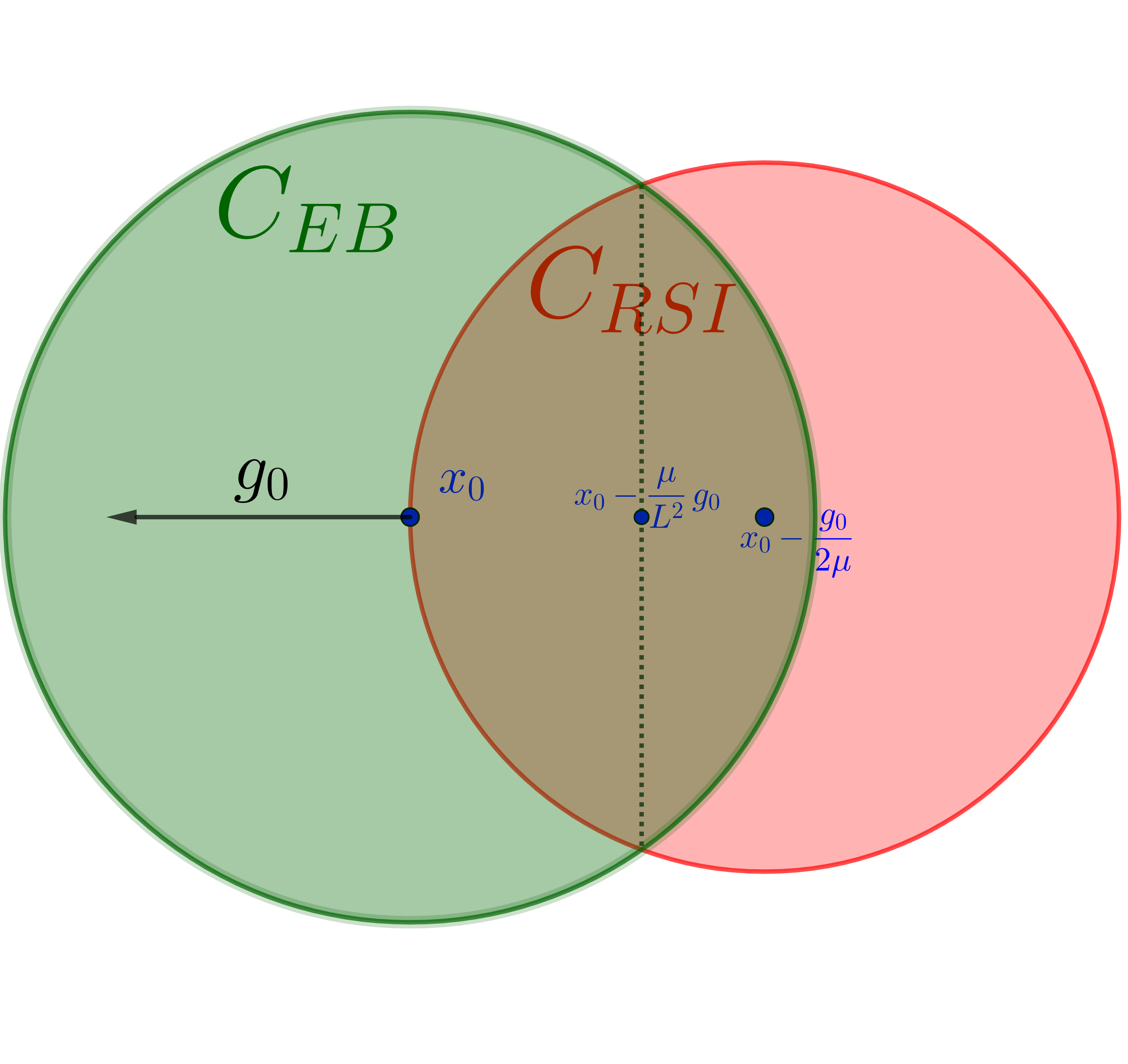}}
\subfigure[$\frac{L}{\mu}=\sqrt{2}$]{\label{fig:2}\includegraphics[width=39mm]{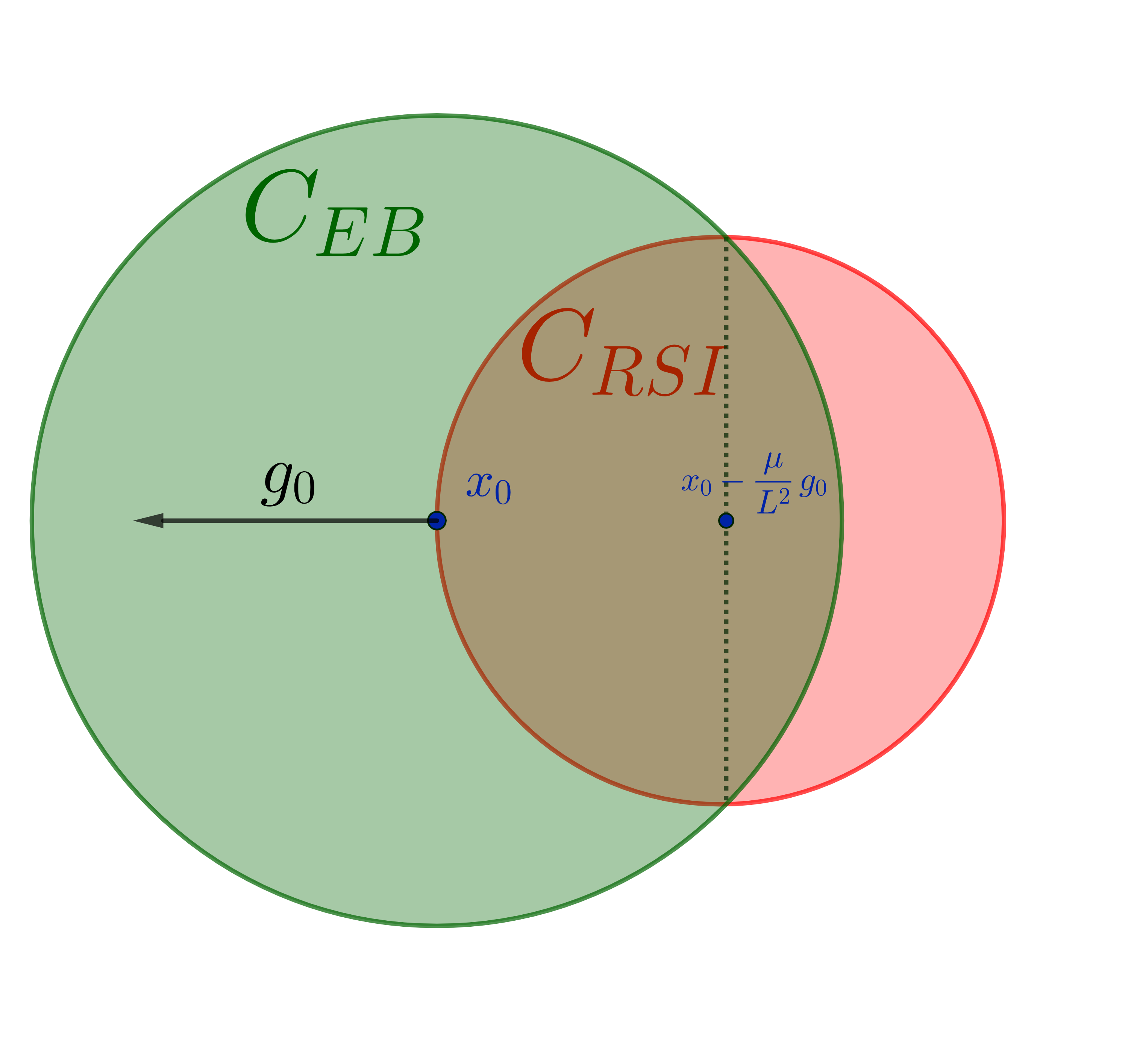}}
\subfigure[$\frac{L}{\mu}<\sqrt{2}$]{\label{fig:3}\includegraphics[width=39mm]{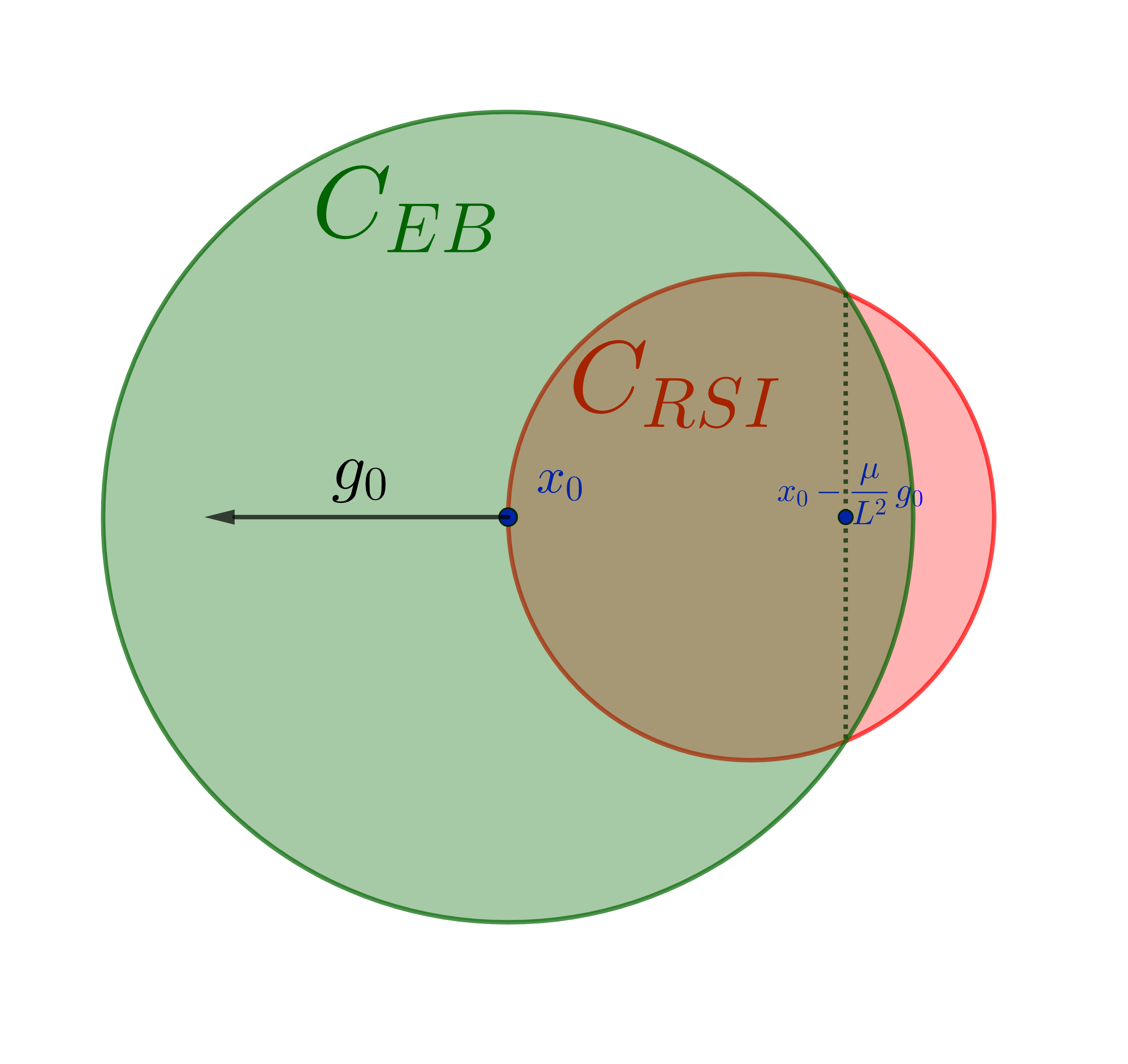}}
\caption{2D representations of the possible positions of $x^*$ given $x_0$ and $g_0$, for three possible values of $\frac{L}{\mu}$. $x^*$ must be in $C_{RSI}$ (red) but not in $C_{EB}$ (green). When $\frac{L}{\mu}> \sqrt{2}$ (left), doing a larger step to reach the center of $C_{RSI}$ will minimize worst-case sub-optimality.}
\vspace{-2mm}
\label{fig:geom}
\end{figure}

Finally, we propose a geometric interpretation of the threshold $\frac{L}{\mu}=\sqrt{2}$. Given $x_0$ and $g_0$, $\RSI^-(\mu)$ requires $x^*$ to be within the circle of center $x_0-\frac{g_0}{2\mu}$ and radius $\frac{\le\|g_0\ri\|_2}{2\mu}$, while $\EB^+(L)$ requires $x^*$ to \emph{not} be within the circle of center $x_0$ and radius $\frac{\le\|g_0\ri\|_2}{L}$. In Figure \ref{fig:geom} we show these circles for different values of $\frac{L}{\mu}$. When $\frac{L}{\mu}>\sqrt{2}$ (Figure \ref{fig:1}), $x_1=x_0-\frac{\mu}{L^2}g_0$ minimizes $\frac{\le\|x_1-x^*\ri\|_2}{\le\|x_0-x^*\ri\|_2}$ over all possible $x^*$, while $x_1=x_0-\frac{g_0}{2\mu}$ minimizes $\le\|x_1-x^*\ri\|_2$. As $\frac{L}{\mu}$ becomes smaller than $\sqrt{2}$ (Figure \ref{fig:2} and \ref{fig:3}), the same point $x_1=x_0-\frac{\mu}{L^2}g_0$ minimizes both quantities.

\end{subsection}
\begin{subsection}{PEP experiment}

Since we have found necessary and sufficient interpolation conditions in Theorem \ref{interpolation_th}, we can use the PEP framework on $\RSI^-\cap\EB^+$ to confirm our results and derive the worst-case convergence rate of first-order algorithms. In Figure \ref{fig:pep} we show the worst-case linear rate of convergence of Heavy Ball (HB) \citep{polyak_gradient_1963} on $\RSI^-(0.1)\cap\EB^+(1.0)$ for regularly sampled learning rate $\alpha$ and momentum $\beta$ (bright yellow means no linear convergence), generated with PEPit \citep{goujaud2022pepit}. We remind the update rule of HB $$x_{n+1}=x_n-\alpha \nabla f(x_n)+\beta(x_n-x_{n-1})$$ Since gradient descent is a special case of HB where $\beta=0$, we observe as expected that the optimal rate of convergence is achieved for $\beta=0$ and $\alpha=0.1=\frac{\mu}{L^2}$. Moreover, Figure \ref{fig:pep} shows that momentum does not do well on $\RSI^-\cap\EB^+$ but gradient decent benefits from a relative robustness to the tuning of $\alpha$ : we get similar convergence rates for any $\alpha\in[0.05,0.15]$.

\begin{figure}
    \centering
    \includegraphics[scale=0.60]{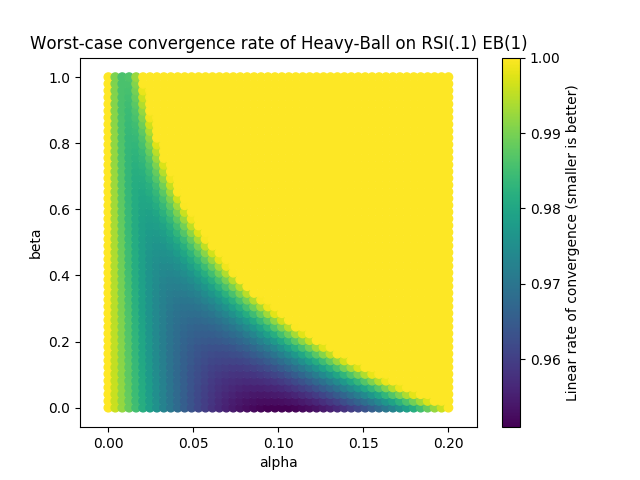}
    \caption{Worst-case linear convergence rate of heavy ball on $\RSI^-(0.1)\cap \EB^+(1)$ depending of its hyperparameters $\alpha$ and $\beta$, as calculated by PEP. The best rate is achieved for $\alpha=0.1$ and $\beta=0$.}
    \label{fig:pep}
    \vspace{-2mm}
\end{figure}

\end{subsection}

\vspace{-2mm}
\section{Conclusion}

    \label{conclusion}

Our main result is to prove that for any $\mu>0$ and $L\geq \mu$, gradient descent is exactly optimal on the class of functions $\RSI^-(\mu)\cap \EB^+(L)$ (by exact optimality, we mean that the convergence guarantees of GD match the lower bound of worst-case performances exactly, without a constant factor of difference). This result confirms the observation in \citep{guille2021study} that optimality is overly sensitive to the choice of assumptions, and should thus be considered with a lot of caution. 

Interestingly, our analysis also identifies two similar notions of optimality, one of which suggests using a larger step size on the first iteration when the function is not particularly well-conditioned to improve worst-case convergence speed (see section \ref{sec:discussion}).

We verified empirically that $RSI^-\cap EB^+$, with respect to the last iterate, are verified on the optimization paths of simple deep neural networks (c.f. Appendix \ref{appendix:rsiebexp}). This suggests that unlike usual alternatives which are known to not be verified on the highly non-convex loss landscapes of neural networks, convergence guarantees on $RSI^-\cap EB^+$ realistically apply to deep learning. On the other hand, the impossibility to accelerate the convergence rate on $RSI^-\cap EB^+$ implies that these assumptions are insufficient to explain the empirical successes of accelerated methods. 

For the scope of this work we have focused on worst-case convergence analysis. While in practice average-case convergence rates are more insightful than their worst-case counterparts, such results are rare due to the necessity of defining a reasonable distribution on the considered class of functions, which is generally unfeasible. While such distribution on $\RSI^-\cap\EB^+$ is equally difficult to define, it should be feasible to define instead a reasonable distribution of the observed gradient for a given sampling point $x$ and nearest minima $x^*$, e.g. an uniform distribution over the (simple) set of possible gradients. Such approach is conceivable with $\RSI^-\cap\EB^+$ only because the interpolation conditions are independent (see remark \ref{rem:indep}), and thus we can easily make sure that any set of gradients sampled from this distribution can be interpolated within the class. While the distribution will necessarily be arbitrary, we believe such analysis could yield very useful insights and $\RSI^-\cap\EB^+$ is a rare opportunity to follow this approach.

Finally, many alternative conditions have been introduced in the literature (see Section \ref{related_work}), for which optimality results are still unknown. The PEP framework is a powerful tool to study conditions for which we can determine sufficient and necessary interpolation conditions, which would improve our understanding of first-order algorithm properties and tuning on a wide variety of objective function classes.

\begin{ack}

The authors would like to thank Leonard Boussioux, for useful discussions and feedback. Ioannis Mitliagkas acknowledges support by an NSERC Discovery grant (RGPIN-2019-06512), a Samsung grant and a Canada CIFAR AI chair.
\end{ack}

\bibliographystyle{abbrvnat}
\bibliography{bib}~\nocite{*}

\begin{thebibliography}{40}
\providecommand{\natexlab}[1]{#1}
\providecommand{\url}[1]{\texttt{#1}}
\expandafter\ifx\csname urlstyle\endcsname\relax
  \providecommand{\doi}[1]{doi: #1}\else
  \providecommand{\doi}{doi: \begingroup \urlstyle{rm}\Url}\fi

\bibitem[Agarwal et~al.(2012)Agarwal, Negahban, and
  Wainwright]{agarwal2011fast}
A.~Agarwal, S.~N. Negahban, and M.~J. Wainwright.
\newblock Fast global convergence of gradient methods for high-dimensional
  statistical recovery.
\newblock \emph{Ann. Statist.}, 40\penalty0 (5):\penalty0 2452--2482, 2012.

\bibitem[Anitescu(2000)]{anitescu2000degenerate}
M.~Anitescu.
\newblock Degenerate nonlinear programming with a quadratic growth condition.
\newblock \emph{SIAM Journal on Optimization}, 10\penalty0 (4):\penalty0
  1116--1135, 2000.

\bibitem[Bolte et~al.(2008)Bolte, Daniilidis, Ley, and
  Mazet]{bolte2008characterizations}
J.~Bolte, A.~Daniilidis, O.~Ley, and L.~Mazet.
\newblock Characterizations of {\l}ojasiewicz inequalities and applications.
\newblock \emph{arXiv:0802.0826}, 2008.

\bibitem[Bonnans and Ioffe(1993)]{bonnans1993second}
J.~F. Bonnans and A.~D. Ioffe.
\newblock \emph{Second-order sufficiency and quadratic growth for non isolated
  minima}.
\newblock PhD thesis, INRIA, 1993.

\bibitem[Bubeck(2015)]{bubeck2015convex}
S.~Bubeck.
\newblock Convex optimization: Algorithms and complexity, 2015.

\bibitem[Drori and Taylor(2019)]{2019}
Y.~Drori and A.~B. Taylor.
\newblock Efficient first-order methods for convex minimization: a constructive
  approach.
\newblock \emph{Mathematical Programming}, 184\penalty0 (1-2):\penalty0
  183–220, Jun 2019.
\newblock ISSN 1436-4646.
\newblock \doi{10.1007/s10107-019-01410-2}.
\newblock URL \url{http://dx.doi.org/10.1007/s10107-019-01410-2}.

\bibitem[Drori and Teboulle(2012)]{drori2012performance}
Y.~Drori and M.~Teboulle.
\newblock Performance of first-order methods for smooth convex minimization: a
  novel approach, 2012.

\bibitem[Gong and Ye(2014)]{gong2014linear}
P.~Gong and J.~Ye.
\newblock Linear convergence of variance-reduced stochastic gradient without
  strong convexity.
\newblock \emph{arXiv:1406.1102}, 2014.

\bibitem[Goujaud et~al.(2022)Goujaud, Moucer, Glineur, Hendrickx, Taylor, and
  Dieuleveut]{goujaud2022pepit}
B.~Goujaud, C.~Moucer, F.~Glineur, J.~Hendrickx, A.~Taylor, and A.~Dieuleveut.
\newblock Pepit: computer-assisted worst-case analyses of first-order
  optimization methods in python, 2022.

\bibitem[Gower(2018)]{gower2018convergence}
R.~M. Gower.
\newblock Convergence theorems for gradient descent.
\newblock \emph{Lecture notes for Statistical Optimization}, 2018.

\bibitem[Gower et~al.(2020)Gower, Sebbouh, and Loizou]{gower2020sgd}
R.~M. Gower, O.~Sebbouh, and N.~Loizou.
\newblock Sgd for structured nonconvex functions: Learning rates, minibatching
  and interpolation.
\newblock \emph{arXiv:2006.1031}, 2020.

\bibitem[Guille-Escuret et~al.(2021)Guille-Escuret, Girotti, Goujaud, and
  Mitliagkas]{guille2021study}
C.~Guille-Escuret, M.~Girotti, B.~Goujaud, and I.~Mitliagkas.
\newblock A study of condition numbers for first-order optimization.
\newblock In \emph{International Conference on Artificial Intelligence and
  Statistics}, pages 1261--1269. PMLR, 2021.

\bibitem[Hanzely et~al.(2018)Hanzely, Richtarik, and
  Xiao]{hanzely2018accelerated}
F.~Hanzely, P.~Richtarik, and L.~Xiao.
\newblock Accelerated {B}regman proximal gradient methods for relatively smooth
  convex optimization.
\newblock Technical Report MSR-TR-2018-22, Microsoft, 2018.

\bibitem[Hardt et~al.(2018)Hardt, Ma, and Recht]{hardt2016gradient}
M.~Hardt, T.~Ma, and B.~Recht.
\newblock Gradient descent learns linear dynamical systems.
\newblock \emph{Journal of Machine Learning Research}, 19, 2018.

\bibitem[Hazan et~al.(2015)Hazan, Levy, and Shalev-Shwartz]{NIPS2015_5718}
E.~Hazan, K.~Levy, and S.~Shalev-Shwartz.
\newblock Beyond convexity: Stochastic quasi-convex optimization.
\newblock In C.~Cortes, N.~D. Lawrence, D.~D. Lee, M.~Sugiyama, and R.~Garnett,
  editors, \emph{Advances in Neural Information Processing Systems 28}, pages
  1594--1602. Curran Associates, Inc., 2015.

\bibitem[He et~al.(2015)He, Zhang, Ren, and
  Sun]{https://doi.org/10.48550/arxiv.1512.03385}
K.~He, X.~Zhang, S.~Ren, and J.~Sun.
\newblock Deep residual learning for image recognition, 2015.
\newblock URL \url{https://arxiv.org/abs/1512.03385}.

\bibitem[Hu and Lessard(2017)]{hu2017dissipativity}
B.~Hu and L.~Lessard.
\newblock Dissipativity theory for nesterov's accelerated method, 2017.

\bibitem[Ioffe(1994)]{ioffe1994sensitivity}
A.~Ioffe.
\newblock On sensitivity analysis of nonlinear programs in banach spaces: the
  approach via composite unconstrained optimization.
\newblock \emph{SIAM Journal on Optimization}, 4\penalty0 (1):\penalty0 1--43,
  1994.

\bibitem[Karimi et~al.(2016)Karimi, Nutini, and
  Schmidt]{DBLP:journals/corr/KarimiNS16}
H.~Karimi, J.~Nutini, and M.~Schmidt.
\newblock Linear convergence of gradient and proximal-gradient methods under
  the polyak-{\l}ojasiewicz condition.
\newblock \emph{CoRR}, abs/1608.04636, 2016.
\newblock URL \url{http://arxiv.org/abs/1608.04636}.

\bibitem[Krizhevsky(2009)]{Krizhevsky09learningmultiple}
A.~Krizhevsky.
\newblock Learning multiple layers of features from tiny images.
\newblock Technical report, 2009.

\bibitem[Kurdyka(1998)]{kurdyka1998gradients}
K.~Kurdyka.
\newblock On gradients of functions definable in o-minimal structures.
\newblock \emph{Annales de l'institut Fourier}, 48:\penalty0 769--783, 1998.

\bibitem[Liu and Wright(2015)]{liu2014asynchronous}
J.~Liu and S.~J. Wright.
\newblock Asynchronous stochastic coordinate descent: Parallelism and
  convergence properties.
\newblock \emph{SIAM J. Optim.}, 25\penalty0 (1):\penalty0 351---376, 2015.

\bibitem[Lu et~al.(2018)Lu, Freund, and Nesterov]{lu2016relativelysmooth}
H.~Lu, R.~M. Freund, and Y.~Nesterov.
\newblock Relatively-smooth convex optimization by first-order methods, and
  applications.
\newblock \emph{SIAM J. Optim.}, 28\penalty0 (1):\penalty0 333--354, 2018.

\bibitem[Luo and Tseng(1993)]{luo_error_1993}
Z.-Q. Luo and P.~Tseng.
\newblock Error bounds and convergence analysis of feasible descent methods: a
  general approach.
\newblock \emph{Annals of Operations Research}, 46\penalty0 (1):\penalty0
  157--178, 1993.

\bibitem[Ma et~al.(2016)Ma, Tappenden, and Tak{\`a}\v{c}]{ma2015linear}
C.~Ma, R.~Tappenden, and M.~Tak{\`a}\v{c}.
\newblock Linear convergence of the randomized feasible descent method under
  the weak strong convexity assumption.
\newblock \emph{Journal of Machine Learning Research}, 17\penalty0
  (228):\penalty0 1--24, 2016.

\bibitem[Necoara et~al.(2016)Necoara, Nesterov, and Glineur]{necoara2016linear}
I.~Necoara, Y.~Nesterov, and F.~Glineur.
\newblock Linear convergence of first order methods for non-strongly convex
  optimization, 2016.

\bibitem[Nesterov(1983)]{Nesterov}
Y.~Nesterov.
\newblock A method of solving a convex programming problem with convergence
  rate $\mathcal{O}\left(\frac{1}{k^2}\right)$.
\newblock \emph{Soviet Mathematics Doklady}, 27:\penalty0 372--376, 1983.

\bibitem[Nesterov(2003)]{nesterov2003introductory}
Y.~Nesterov.
\newblock \emph{Introductory lectures on convex optimization: A basic course},
  volume~87.
\newblock Springer Science \& Business Media, 2003.

\bibitem[Polyak(1963)]{polyak_gradient_1963}
B.~T. Polyak.
\newblock Gradient methods for the minimisation of functionals.
\newblock \emph{USSR Computational Mathematics and Mathematical Physics},
  3\penalty0 (4):\penalty0 864 -- 878, 1963.

\bibitem[Polyak(1987)]{polyak1987introduction}
B.~T. Polyak.
\newblock Introduction to optimization. optimization software.
\newblock \emph{Inc., Publications Division, New York}, 1:\penalty0 32, 1987.

\bibitem[Sch{\"o}pfer(2016)]{Schpfer2016LinearCO}
F.~Sch{\"o}pfer.
\newblock Linear convergence of descent methods for the unconstrained
  minimization of restricted strongly convex functions.
\newblock \emph{SIAM J. Optim.}, 26:\penalty0 1883--1911, 2016.

\bibitem[Taylor(2017)]{Taylor2017ConvexIA}
A.~B. Taylor.
\newblock Convex interpolation and performance estimation of first-order
  methods for convex optimization, 2017.

\bibitem[Taylor et~al.(2016)Taylor, Hendrickx, and Glineur]{taylor2016smooth}
A.~B. Taylor, J.~M. Hendrickx, and F.~Glineur.
\newblock Smooth strongly convex interpolation and exact worst-case performance
  of first-order methods, 2016.

\bibitem[Taylor et~al.(2017)Taylor, Hendrickx, and Glineur]{2017}
A.~B. Taylor, J.~M. Hendrickx, and F.~Glineur.
\newblock Exact worst-case performance of first-order methods for composite
  convex optimization.
\newblock \emph{SIAM Journal on Optimization}, 27\penalty0 (3):\penalty0
  1283–1313, Jan 2017.
\newblock ISSN 1095-7189.
\newblock \doi{10.1137/16m108104x}.
\newblock URL \url{http://dx.doi.org/10.1137/16M108104X}.

\bibitem[Taylor et~al.(2020)Taylor, Hendrickx, and Glineur]{taylor2020exact}
A.~B. Taylor, J.~M. Hendrickx, and F.~Glineur.
\newblock Exact worst-case convergence rates of the proximal gradient method
  for composite convex minimization, 2020.

\bibitem[Yi et~al.(2019)Yi, Zhang, Yang, Johansson, and
  Chai]{yi2019exponential}
X.~Yi, S.~Zhang, T.~Yang, K.~H. Johansson, and T.~Chai.
\newblock Exponential convergence for distributed smooth optimization under the
  restricted secant inequality condition, 2019.

\bibitem[Yuan et~al.(2016)Yuan, Ling, and Yin]{Yuan2016OnTC}
K.~Yuan, Q.~Ling, and W.~Yin.
\newblock On the convergence of decentralized gradient descent.
\newblock \emph{SIAM J. Optim.}, 26:\penalty0 1835--1854, 2016.

\bibitem[Zhang(2017)]{zhang2017restricted}
H.~Zhang.
\newblock The restricted strong convexity revisited: analysis of equivalence to
  error bound and quadratic growth.
\newblock \emph{Optimization Letters}, 11\penalty0 (4):\penalty0 817--833,
  2017.

\bibitem[Zhang and Yin(2013)]{zhang2013gradient}
H.~Zhang and W.~Yin.
\newblock Gradient methods for convex minimization: better rates under weaker
  conditions.
\newblock Cam report, UCLA, 2013.

\bibitem[Zhou et~al.(2019)Zhou, Liang, and Shen]{zhou2019simple}
Y.~Zhou, Y.~Liang, and L.~Shen.
\newblock A simple convergence analysis of bregman proximal gradient algorithm.
\newblock \emph{Computational Optimization and Applications}, 73\penalty0
  (3):\penalty0 903--912, 2019.

\end{thebibliography}

\newpage
\appendix

\section{\texorpdfstring{Optimization paths of neural networks are in $RSI^-\cap EB^+$}{Optimization paths of neural networks are in RSI- and EB+}}
    \label{appendix:rsiebexp}

    In this section we present a short experiment motivating $RSI^-\cap EB^+$ as useful assumptions for the study of optimization of neural networks. The goal is to estimate whether the gradients seen by training neural networks are interpolable by a function $f\in RSI^-\cap EB^+$. We propose the following process :
\begin{enumerate}
    \item Set random seed S
    \item Train a ResNet18 on CIFAR10 until convergence, and save the last iterate $x^*$
    \item Reset random seed to S
    \item Train a ResNet18 on CIFAR10, and at each iteration $i$, sample batch $B$ and measure $RSI_i=\frac{\langle \nabla f_B(x_i)\mid x_i-x^*\rangle}{\|x_i-x^*\|_2^2}$ and $EB_i=\frac{\|\nabla f_B(x_i)\|_2}{\|x_i-x^*\|_2}$
\end{enumerate}
Where $\nabla f_B(x)$ is the gradient of the loss function on minibatch $B$, $x_i$ is the value of the weight at iteration $i$, and $x^*$ is the value of the last iterate measured at step 2. 

Due to resetting the seed to a same value, the two training runs will be identical. We consider the last iterate $x^*$ to approximate a local minima, and $RSI_i$ and $EB_i$ will indicate whether the gradients seen during optimization are compatible with $RSI^-\cap EB^+$ with respect to that minima.

Even if the full-batch objective function that we intend to optimize is in $RSI^-\cap EB^+$, it is possible for $RSI_i$ to be negative due to the variance w.r.t the sampling of the minibatch $B$. However, we observe empirically that despite this, $RSI_i$ is lower bounded by a strictly positive value for every single iteration, without exception. This behavior is consistent across optimization algorithms (LARS, SGD without momentum, SGD with momentum, and ADAM) and initializations. For simplicity, we present here the results using SGD without momentum, with learning rate $0.1$ and batch size $|B|=1000$, and train for $360$ epochs. 

\begin{figure}[H]

\centering
\includegraphics[width=0.6\textwidth]{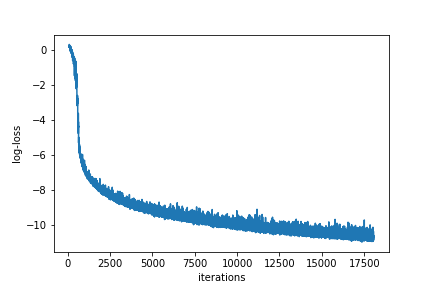}\hfill

\caption{log-loss throughout training.}
\label{fig:appendix:loss}

\end{figure}

\begin{figure}[H]

\centering
\includegraphics[width=.5\textwidth]{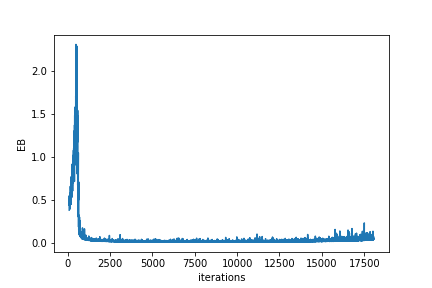}\hfill
\includegraphics[width=.5\textwidth]{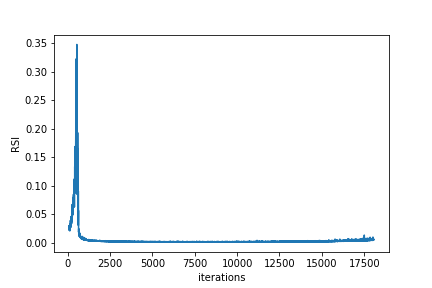}\hfill

\caption{$RSI_i$ (right) and $EB_i$ (left) throughout training from epoch $0$ to $360$.}
\label{fig:appendix:rsieb}

\end{figure}

\begin{figure}[H]

\centering
\includegraphics[width=.5\textwidth]{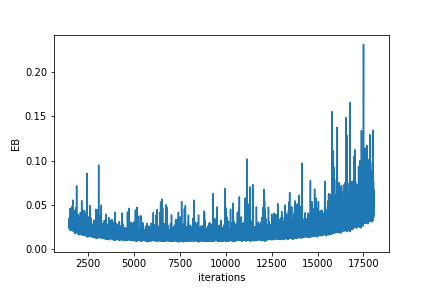}\hfill
\includegraphics[width=.5\textwidth]{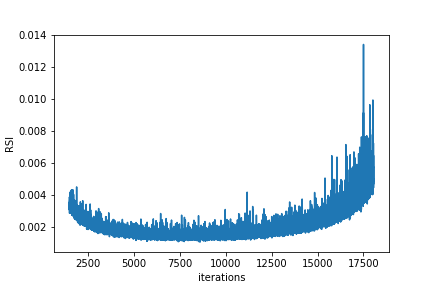}\hfill

\caption{$RSI_i$ (right) and $EB_i$ (left) throughout training from epoch $30$ to $360$.}
\label{fig:appendix:rsieb_30}

\end{figure}

\textbf{Results:} we report the log of training loss in Figure \ref{fig:appendix:loss} and the measured $RSI_i$ and $EB_i$ in Figure \ref{fig:appendix:rsieb}. Moreover, in order to better observe the behavior of $RSI_i$ and $EB_i$ outside of the initial peak, we report in Figure \ref{fig:appendix:rsieb_30} the measured $RSI_i$ and $EB_i$ starting at epoch $30$. We observe that $EB_i$ is upper bounded by $L=2.303$ and $RSI_i$ is lower bounded by $\mu=0.0010$, with a resulting condition number $\kappa=\frac{L}{\mu}=2196.0$. Both have a significant peak at the beginning of training, justifying the popular use of learning rate warm-ups. When measuring bounds after epoch $30$, we obtain $L=0.2308$ resulting in a condition number $\kappa=220.1$. 

Surprisingly, despite the variance induced by minibatch sampling, the observed $RSI_i$ are all lower bounded by $\mu=0.0010>0$. In particular, due to the necessary and sufficient conditions of $RSI^-\cap EB^+$ (See Section \ref{interpolation}), it is guaranteed that there exists a function $f\in RSI^-\cap EB^+$ which exactly interpolates the gradients seen by the optimizer. And therefore, the convergence guarantees of $RSI^-\cap EB^+$ naturally apply to the optimization of neural networks in this setting. 

Note that we do not claim that the objective function is in $RSI^-\cap EB^+$, which seems unlikely, but that the iterates explored by first-order algorithms are interpolable by functions in $RSI^-\cap EB^+$, including when sampling only part of the objective function through minibatches. This result strongly motivates the study of $RSI^-\cap EB^+$ as its guarantees apply to the optimization of neural network under assumptions empirically verified (at least in this simple setting).


\section{Proof of Theorem \ref{interpolation_th}}

    \label{proof_interpolation}
We start with two simple lemmas
\begin{lem} \label{lem:proof_thm1_lemma_project_inequality}
Let $X^*$ be a closed convex set, and $x^* \in X^*$ be the orthogonal projection of $x$ onto $X^*$. Then for any $y \in X^*$, 
\begin{align}
    \braket{x^* - y}{x - x^*} \geq 0
\end{align}
\end{lem}
\textbf{Proof.} Let $y \in X^*$.

For $\theta\in[0,1]$, \\
let $h(\theta) = \norm{x - \lrpar{(1-\theta) x^* + \theta y}}^2 = \le\|x-x^*\ri\|_2^2+2\theta\braket{x - x^*}{x^* - y} +\theta^2\le\|x^*-y\ri\|_2^2$.

$h$ is differentiable and
\begin{align}
    \label{proof_eq_deriv1}
    h'(\theta) &= 2 \braket{x - x^*}{x^* - y} + 2\theta\le\|x^*-y\ri\|_2^2
\end{align}

Since $x^*$ is the orthogonal projection of $x$ onto $X^*$ and $\forall\theta\in[0,1], (1-\theta)x^*+\theta y\in X^*$, we have $\forall\theta\in[0,1], h(\theta)\geq h(0)$, and thus $h'(0)\geq 0$. This concludes the proof of the Lemma thanks to \eqref{proof_eq_deriv1}.

$\hfill\blacksquare$

\begin{lem}
If $x$ and $x_i$ are two points with respective orthogonal projections $x^*$ and $x_i^*$ on a closed convex set, then 
\begin{align}
    \le\|x-x^*-(x_i-x_i^*)\ri\|_2\leq 2\le\|x-x_i\ri\|_2
\end{align}
\label{lemma_interp}
\end{lem}
\textbf{Proof.} As the case $x^* = x_i^*$ is trivial, we may assume that $x^* \neq x_i^*$. 

Using lemma \ref{lem:proof_thm1_lemma_project_inequality} twice, we get
\begin{align}
    0 &\leq \braket{x - x^*}{x^* - x_i^*} \\
    0 &\leq \braket{x_i - x_i^*}{x_i^* - x^*} = \braket{x_i^* - x_i}{x^* - x_i^*}
\end{align}
Adding the two inequalities, we get that
\begin{align}
    0 \leq \braket{x - x^* - x_i + x_i^*}{x^* - x_i^*} &= \braket{x - x_i}{x^* - x_i^*} - \norm{x^* - x_i^*}^2 \nonumber \\
    &\leq \norm{x - x_i} \norm{x^* - x_i^*} - \norm{x^* - x_i^*}^2
\end{align}
Since $x^* \neq x_i^*$, we obtain
\begin{align}
    \norm{x^* - x_i^*} \leq \norm{x - x_i}
\end{align}
And thus
\begin{equation}
    \norm{x - x^* - (x_i - x_i^*)} \leq \norm{x-x_i}+\norm{x^*-x_i^*}\leq 2\norm{x - x_i}
\end{equation}

$\hfill\blacksquare$

We now move on to the proof of Theorem \ref{interpolation_th}, that is :

Let $(x_i, g_i)_{i\leq n} \in \le(\R^d\times\R^d\ri)^{n+1}$, such that the $x_i$ are separate points.

Then, $\forall \mu, L>0$:
$$\exists f\in \RSI^-(\mu)\cap \EB^+(L),\: s.t.\: \forall i, \: \nabla f(x_i)=g_i$$
$$\Updownarrow$$
$$\exists X^*\subseteq \mathbb{R}^d\ convex,\; s.t.\; \forall i,$$
\begin{equation}
    \begin{aligned}
        \le\|g_i\ri\|_2 \leq L\le\|x_i-x_i^*\ri\|_2\quad \text{and}\quad \left\langle g_i \mid x_i-x_i^* \right\rangle \geq \mu\le\|x_i-x_i^*\ri\|_2^2 
    \end{aligned}
\end{equation}

Where $x_i^*$ is the orthogonal projection of $x_i$ onto $X^*$.

\textbf{Proof.} The direct implication is trivial since the second property is simply the application of $\RSI^-$ and $\EB^+$ in each $x_i$. Let us now assume that \eqref{interp} is verified. 

First, let us note that if $L=\mu$, then we have $\forall i, g_i = \mu(x_i-x_i^*)$ and thus we can easily interpolate the $(x_i,g_i)$ using $f(x)=\frac{\mu}{2}\le\|x-x^*\ri\|_2^2$. We now assume $L>\mu$. 


If there is only one pair $(x_i,g_i)$, then we can simply use $f(x)=\le\langle g_i\mid x-x_i\ri\rangle + \frac{\mu+L}{4}\le\|x-x_i\ri\|_2^2$ to interpolate $(x_i,g_i)$. $f$ is then $\mu$-strongly convex and $L$-smooth so it is also in $RSI^-(\mu)$ and $EB^+(L)$. Let us now assume there are at least two pairs $(x_i,g_i)_i$. Let

\begin{equation}
    \eps_0 = \frac{1}{2}\min_{i\neq j}(\le\|x_i-x_j\ri\|_2)>0    
\end{equation}

By construction,

\begin{equation}
\label{unique_i}
\forall x\in\mathbb{R}^d, \le(\exists i, \le\|x-x_i\ri\|_2< \eps_0\ri) \Rightarrow \forall j\neq i, \le\|x-x_j\ri\|_2\geq \eps_0    
\end{equation}

Moreover, if $\forall i, x_i\in X^*$, we can simply take $f(x)=\frac{\mu}{2}\le\|x-x^*\ri\|_2^2$. Otherwise, let $\mathcal{I}=\le\{i\mid x_i\neq x_i^*\ri\}$, and let $\eps_1 = \frac{1}{2}\min_{i\in\mathcal{I}}(\le\|x_i-x_i^*\ri\|_2)>0$\\

Let $\eps<\min(\eps_0,\eps_1)$ and $0<\beta<\frac{1}{2}$. We introduce the function $\lambda_{\eps,\beta}$ from $\le[0,\eps\ri]$ to $\le[0,1\ri]$ defined by :

\begin{equation}
    \lambda_{\eps,\beta}(u) = \frac{1+\cos\le(\pi\frac{u^\beta}{\eps^\beta}\ri)}{2}
\end{equation}

We finally introduce our interpolation function :

\begin{gather}
f_{\eps,\beta}(x) = \begin{cases}
\frac{\mu+L}{4}\le\|x-x^*\ri\|_2^2 & if\; \forall i, \le\|x-x_i\ri\|_2 \geq \eps\\
\frac{\mu+L}{4}\le\|x-x^*\ri\|_2^2 + \lambda_{\eps,\beta}\le(\le\|x-x_i\ri\|_2\ri)\le\langle g_i-\frac{\mu+L}{2}(x_i-x_i^*)\mid x-x_i\ri\rangle & if\; \exists i, \le\|x-x_i\ri\|_2 < \eps
\end{cases}
\label{f_def}
\end{gather}

First let us note that $f_{\eps,\beta}$ is properly defined : as stated in \eqref{unique_i}, there may be at most one $i$ such that $\le\|x-x_i\ri\|_2<\eps$. Moreover, $f_{\eps,\beta}$ is continuous because $\lambda_{\eps,\beta}(\eps)=0$.

Since $\lambda_{\eps,\beta}(\eps)=0$ and $\lambda_{\eps,\beta}'(\eps)=0$, we can easily verify that for $x$ such that $\le\|x-x_i\ri\|_2=\eps$, $f_{\eps,\beta}$ is differentiable in $x$ with $\nabla f_{\eps,\beta}(x)=\frac{\mu+L}{2}(x-x^*)$. Thus $f_{\eps,\beta}$ is differentiable on $\mathbb{R}^d$. 
For any $x\in\mathbb{R}^d$ such that $\forall i, \le\|x-x_i\ri\|_2\geq \eps$, we have $\nabla f_{\eps,\beta}(x)=\frac{\mu+L}{2}(x-x^*)$ and thus trivially 
\begin{equation}
\label{eq6}
   \le\langle\nabla f_{\eps,\beta}(x)\mid x-x^*\ri\rangle = \frac{\mu+L}{2}\le\|x-x^*\ri\|_2^2 \geq \mu\le\|x-x^*\ri\|_2^2
\end{equation}
\begin{equation}
\label{eq7}
   \le\|\nabla f_{\eps,\beta}(x)\ri\| = \frac{\mu+L}{2}\le\|x-x^*\ri\|_2 \leq L\le\|x-x^*\ri\|_2
\end{equation}


Let us now assume there is $i$ such that $\le\|x-x_i\ri\|_2<\eps$. If $x_i=x_i^*$, then $g_i=0$ and $\nabla f(x) = \frac{\mu+L}{4}(x-x^*)$ and equations \eqref{eq6} and \eqref{eq7} are respected as well. Otherwise, we have $\le\|x-x^*\ri\|_2\geq \min_{i\in\mathcal{I}}(\le\|x_i-x_i^*\ri\|_2)-\eps = \eps_1-\eps > 0$ \\
We then have, for $x\neq x_i$:

\begin{equation}
\label{gradf}
\begin{aligned}
\nabla f_{\eps,\beta}(x) & = \frac{\mu+L}{2}(x-x^*) + \lambda_{\eps,\beta}\le(\le\|x-x_i\ri\|_2\ri)\le(g_i-\frac{\mu+L}{2}\le(x_i-x_i^*\ri)\ri)\\
                       & \quad + \lambda_{\eps,\beta}'\le(\le\|x-x_i\ri\|_2\ri)\frac{x-x_i}{\le\|x-x_i\ri\|_2}\le\langle g_i-\frac{\mu+L}{2}(x_i-x_i^*)\mid x-x_i \ri\rangle \\
                       & = \le(1 - \lambda_{\eps,\beta}\le(\le\|x-x_i\ri\|_2\ri)\ri)\frac{\mu+L}{2}(x-x^*) + \lambda_{\eps,\beta}(\le\|x-x_i\ri\|_2)g_i \\
                       & \quad + \lambda_{\eps,\beta}(\le\|x-x_i\ri\|_2)\frac{\mu+L}{2}\le(x-x^*-(x_i-x_i^*)\ri) \\
                       & \quad - \frac{\pi}{2}\frac{\le\|x-x_i\ri\|_2^\beta}{\eps^\beta}\frac{\beta(x-x_i)}{\le\|x-x_i\ri\|_2^2}sin(\pi\frac{\le\|x-x_i\ri\|_2^\beta}{\eps^\beta})\le\langle g_i-\frac{\mu+L}{2}(x_i-x_i^*)\mid x-x_i \ri\rangle
\end{aligned}
\end{equation}

Since $X^*$ is a convex set (and closed by continuity of $f_{\eps,\beta}$), we have from Lemma \ref{lemma_interp} $\le\|(x-x^*-(x_i-x_i^*))\ri\|_2\leq 2\le\|x-x_i\ri\|_2$. 

To simplify notations, let us note $u=\le\|x-x_i\ri\|_2$, $\lambda = \lambda_{\eps,\beta}(u)$, and

$r=-\frac{\pi}{2}\frac{u^\beta}{\eps^\beta}\frac{\beta(x-x_i)}{u^2}sin(\pi\frac{u^\beta}{\eps^\beta})\le\langle g_i-\frac{\mu+L}{2}(x_i-x_i^*)\mid x-x_i \ri\rangle$.

We first want to upper bound $\le\|\nabla f_{\eps,\beta}\ri\|_2$ using \eqref{gradf}: 

\begin{equation}
\label{eq9}
\begin{aligned}
\le\|\nabla f_{\eps,\beta}(x)\ri\|_2 & \leq \le(1-\lambda\ri)\frac{\mu+L}{2}\le\|x-x^*\ri\|_2 + \lambda\le\|g_i\ri\|_2 + (\mu+L)\lambda u + \le\|r\ri\|_2\\
                                     & \leq \le(1-\lambda\ri)\frac{\mu+L}{2}\le\|x-x^*\ri\|_2 + \lambda L\le(\le\|x-x^*\ri\|_2+\le\|x_i-x_i^*\ri\|-\le\|x-x^*\ri\|\ri) + (\mu+L)\lambda u + \le\|r\ri\|_2\\
                                     & \leq L\le\|x-x^*\ri\|_2-(1-\lambda)\frac{L-\mu}{2}\le\|x-x^*\ri\|_2+(\mu+3L)\lambda u + \le\|r\ri\|_2
\end{aligned}
\end{equation}

Moreover, the $3$rd order remainder of the Taylor expansion of $cos(\pi\frac{u^\beta}{\eps^\beta})$ is $\frac{cos(c)}{4!}(\pi^4\frac{u^{4\beta}}{\eps^{4\beta}})$ for some $c$ in $[0,\pi\frac{u}{\eps}]$ and by upper bounding it we get $cos(\pi\frac{u^\beta}{\eps^\beta})\leq 1-\frac{\pi^2 u^{2\beta}}{2\eps^{2\beta}}+\frac{\pi^4 u^{4\beta}}{24 \eps^{4\beta}}$ and thus, for $\frac{u}{\eps}\leq 1$,
\begin{equation}
\label{eq10}
    \begin{aligned}
    -(1-\lambda)\frac{L-\mu}{2}\le\|x-x^*\ri\|_2 & \leq \le(-\frac{\pi^2 u^{2\beta}}{4\eps^{2\beta}}+\frac{\pi^4 u^{4\beta}}{48\eps^{4\beta}}\ri)\frac{L-\mu}{2}\le\|x-x^*\ri\|_2\\
                                                 & \leq - C_0 \frac{u^{2\beta}}{\eps^{2\beta}}
    \end{aligned}
\end{equation}

with $C_0=\le(\frac{\pi^2}{4}-\frac{\pi^4}{48}\ri)\frac{L-\mu}{2}\eps_1 > 0$\\
Furthermore, since $2\beta<1$ and $\frac{u}{\eps}\leq 1$, we can also bound 
\begin{equation}
\label{eq11}
    \begin{aligned}
        (\mu+3L)\lambda u \leq \eps(\mu+3L)\frac{u}{\eps}\leq \eps C_1 \frac{u^{2\beta}}{\eps^{2\beta}}
    \end{aligned}
\end{equation}

with $C_1 = \mu+3L>0$

Finally, we can bound the last term using $sin(x)\leq |x|$

\begin{equation}
\label{eq12}
    \begin{aligned}
    \le\|r\ri\|_2 & \leq \beta\frac{\pi^2}{2}\frac{3L+\mu}{2}\le\|x_i-x_i^*\ri\|_2\frac{u^{2\beta}}{\eps^{2\beta}}\\
                  & \leq \beta C_2 \frac{u^{2\beta}}{\eps^{2\beta}}
    \end{aligned}
\end{equation}

with $C_2=\frac{\pi^2}{2}\frac{3L+\mu}{2}\max_i(\le\|x_i-x_i^*\ri\|_2)$

Finally, by choosing $\eps\leq \frac{C_0}{2C_1}$ and $\beta\leq \frac{C_0}{2C_2}$, and plugging \eqref{eq10}, \eqref{eq11}, \eqref{eq12} into \eqref{eq9}, we get :

$$\le\|\nabla f_{\eps,\beta}(x)\ri\|_2 \leq L\le\|x-x^*\ri\|_2$$

It only remains now to adequately lower bound $\le\langle \nabla f_{\eps,\beta}(x)\mid x-x^*\ri\rangle$. We use the same method as before and keep the notations :

\begin{equation}
\label{eq13}
    \begin{aligned}
     \le\langle \nabla f_{\eps,\beta}(x)\mid x-x^*\ri\rangle & = (1-\lambda)\frac{\mu+L}{2}\le\|x-x^*\ri\|_2^2 + \lambda\le\langle g_i\mid x-x^*\ri\rangle \\
     & + \lambda \frac{\mu+L}{2}\le\langle x-x^*-(x_i-x_i^*)\mid x-x^*\ri\rangle + \le\langle r\mid x-x^*\ri\rangle\\
     & \geq (1-\lambda)\frac{\mu+L}{2}\le\|x-x^*\ri\|_2^2 + \lambda\mu\le\|x_i-x_i^*\ri\|_2^2+\lambda\le\langle g_i\mid x-x^*-(x_i-x_i^*)\ri\rangle\\
     & - \lambda(\mu+L)u\le\|x-x^*\ri\|_2 - \le\|r\ri\|_2\le\|x-x^*\ri\|_2\\
     & \geq (1-\lambda)\frac{\mu+L}{2}\le\|x-x^*\ri\|_2^2-(1-\lambda)\mu\le\|x-x^*\ri\|_2^2 + \mu\le\|x-x^*\ri\|_2^2 \\
     & + \lambda\mu\le(\le\|x_i-x_i^*\ri\|_2^2-\le\|x-x^*\ri\|_2^2\ri) - 2u\lambda\le\|g_i\ri\|_2 \\
     & - \lambda(\mu+L)u\le\|x-x^*\ri\|_2 - \le\|r\ri\|_2\le\|x-x^*\ri\|_2\\
     & \geq \mu\le\|x-x^*\ri\|_2^2 + C_0\frac{u^{2\beta}}{\eps^{2\beta}} - 2\mu u (2\le\|x_i-x_i^*\ri\|_2+u) \\
     & - 2\eps\lambda L\le\|x_i-x_i^*\ri\|_2\frac{u}{\eps}-\eps\lambda(\mu+L)\le\|x-x^*\ri\|_2\frac{u}{\eps} - \beta C_2 \le\|x-x^*\ri\|_2\frac{u^{2\beta}}{\eps^{2\beta}}\\
     & \geq \mu\le\|x-x^*\ri\|_2^2 + (C_0 - \eps M_0 - \beta M_1)\frac{u^{2\beta}}{\eps^{2\beta}}
    \end{aligned}
\end{equation}

With $M_0=4(\mu+L)\max_i(\le\|x_i-x_i^*\ri\|_2)+(L+3\mu)\eps_1> 0$ \\
and $M_1=C_2(\max_i(\le\|x_i-x_i^*\ri\|_2)+\eps_1)> 0$\\

Therefore by taking $\eps\leq \frac{C_0}{2M_0}$ and $\beta\leq \frac{C_0}{2M_1}$, we guarantee 

$$\le\langle \nabla f_{\eps,\beta}(x)\mid x-x^*\ri\rangle \geq \mu \le\|x-x^*\ri\|_2^2$$

Finally, for any $i$, we have :

\begin{equation}
    \begin{aligned}
    \le\|x-x_i\ri\|_2\lambda'_{\eps,\beta}(\le\|x-x_i\ri\|_2)=-\frac{\pi}{2}\frac{\le\|x-x_i\ri\|_2^\beta}{\eps^{\beta}}\frac{\beta(x-x_i)}{\le\|x-x_i\ri\|_2}\sin(\pi\frac{\le\|x-x_i\ri\|_2^\beta}{\eps^\beta})
    \end{aligned}
\end{equation}

Which goes to $0$ as $x$ tends to $x_i$. Therefore using the definition of $f_{\eps,\beta}$ in \eqref{f_def} and the fact that $\le\langle g_i-\frac{\mu+L}{2}(x_i-x_i^*)\mid x-x_i\ri\rangle$ is linear in $x-x_i$, we can conclude $\nabla f_{\eps,\beta}(x_i)=\frac{\mu+L}{2}(x_i-x_i^*)+\lambda_{\eps,\beta}(0)(g_i-\frac{\mu+L}{2}(x_i-x_i^*))=g_i$.

We thus have proven that for sufficiently small $\eps$ and $\beta$, $\forall i, \nabla f_{\eps,\beta}(x_i)=g_i$, that for all $x$, $\le\langle \nabla f_{\eps,\beta}(x)\mid x-x^*\ri\rangle\geq\mu\le\|x-x^*\ri\|_2^2$ and $\le\|\nabla f_{\eps,\beta}(x)\ri\|_2\leq L\le\|x-x^*\ri\|_2$. Therefore by definition, $f_{\eps,\beta}$ is in $\RSI^-(\mu)\cap \EB^+(L)$ and interpolates the $x_i, g_i$.

Which concludes the proof.

$\hfill\blacksquare$

\section{Proof of Lemma \ref{lemma_lower}}

    \label{lowerbound_proof}
Let $\mu>0$ and $L>\mu$. Let $\alpha_0\in\le[\frac{\mu}{L^2},\max\le(\frac{\mu}{L^2},\frac{1}{2\mu}\ri)\ri]$. For any first-order optimization algorithm $\mathcal{A}$ and starting point $x_0\in\R^d$, there exists $(g_i)_{i\leq (d-2)}\in\mathbb{R}^{d}$, $(f_i)_{i\leq (d-2)}\in\R$ and  $\mathcal{S}_{d-2}\subseteq\mathcal{S}_{d-1}\subseteq\dots\subseteq\mathcal{S}_0\subseteq\R^{d}$ such that:
\begin{enumerate}
    \item $\forall i\leq d-2$, there exists a $(d-i-1)$-dimensional affine space $\mathcal{H}_i$ containing $\mathcal{S}_i$ and in which $\mathcal{S}_i$ is a $(d-i-2)$-sphere of radius $r_i=\sqrt{\frac{\alpha_0}{\mu}-\alpha_0^2}\le\|g_0\ri\|_2\le(1-\frac{\mu^2}{L^2}\ri)^{\frac{i}{2}}$ and center $c_i\in\mathcal{H}_i$.
    \item Let $(x_i)_i$ the iterates generated by $\mathcal{A}$ starting from $x_0$ and reading gradients $(g_i)_i$ and function values $(f_i)_i$, then for any $i\leq d-2$ and any $x\in\mathcal{S}_i$, there exists a function $f$ in $\RSI^-(\mu)\cap \EB^+(L)$ minimized by $\{x\}$ that interpolates $(x_j,f_j,g_j)_{j\leq i}$.
\end{enumerate}
\textbf{Proof.} For any first-order optimization algorithm $\mathcal{A}$ and starting point $x_0$, we are going to construct by induction the sequences $(g_i)_{i}$, $(f_i)_i$ and  $(\mathcal{S}_i)_i$.

\textbf{Initialisation:} Let $g_0\in\R^d\setminus\{0\}$. We can take any non-zero gradient as initialisation. Let $c_0=x_0-\alpha_0 g_0$, $f_0=\frac{\mu+L}{4\mu}\alpha_0\le\|g_0\ri\|_2^2$, and 
$$\mathcal{H}_0=\{x\in\R^d\mid \le\langle x-c_0\mid g_0\ri\rangle = 0\}$$
$\mathcal{H}_0$ is an hyperplane with dimension $d-1$, and we finally introduce 
$$\mathcal{S}_0 = \le\{x\in \mathcal{H}_0\mid \le\|x-c_0\ri\|_2=\sqrt{\frac{\alpha_0}{\mu}-\alpha_0^2}\le\|g_0\ri\|_2\ri\}$$
By construction, $c_0\in\mathcal{H}$ and $\mathcal{S}_0$ is the $(d-2)$-sphere in $\mathcal{H}_0$ of center $c_0$ and radius $r_0=\sqrt{\frac{\alpha_0}{\mu}-\alpha_0^2}\le\|g_0\ri\|_2$.
Moreover, let $x^*\in\mathcal{S}_0$. We have

$$\le\|x^*-x_0\ri\|_2^2=\le\|x^*-c_0+c_0-x_0\ri\|_2^2=r_0^2+\alpha_0^2\le\|g_0\ri\|_2^2=\frac{\alpha_0\le\|g_0\ri\|_2^2}{\mu}$$

And thus $f_0=\frac{\mu+L}{4}\le\|x^*-x_0\ri\|_2^2$. We also have :
$$\le\langle g_0\mid x_0-x^*\ri\rangle = \le\langle g_0\mid x_0-c_0\ri\rangle + \le\langle g_0\mid c_0-x^*\ri\rangle = \alpha_0\le\|g_0\ri\|_2^2+0=\mu\le\|x_0-x^*\ri\|_2^2$$
Finally, since $\alpha_0\geq\frac{\mu}{L^2}$, $$\le\|g_0\ri\|_2^2=\frac{\mu}{\alpha_0}\le\|x_0-x^*\ri\|_2^2\leq L^2\le\|x_0-x^*\ri\|_2^2$$.
Therefore all the sufficient conditions of Corollary \ref{interp_cor} are verified, and there exists $f\in RSI^-(\mu)\cap EB^+(L)$ which is minimized by $\le\{x^*\ri\}$ and interpolates $(x_0,f_0,g_0)$. This concludes the initialization.

\textbf{Induction:} Let us assume the existence of such $(f_j)_j$, $(g_j)_j$ and $(\mathcal{S}_j)_j$ up to step $i\leq d-3$. Let $x_{i+1}$ be the iterate given by $\mathcal{A}$ after reading iterates $(x_j)_j$, function values $(f_j)_j$, and gradients $(g_j)_j$, let $\mathcal{H}_i$ the $(d-i-1)$-dimensional affine space in which $\mathcal{S}_i$ is a sphere, and let $c_i\in\mathcal{H}_i$ the center of the sphere $\mathcal{S}_i$.

If there exists $j\leq i$ such that $x_{i+1}=x_j$, then we simply return $g_{i+1}=g_j$ and $f_{i+1}=f_j$. We can take as $S_{i+1}$ any $(d-i-2)$-dimensional sphere of radius $r_{i+1}$ included in $S_i$, $\mathcal{H}_{i+1}$ its supporting affine space and $c_{i+1}$ its center. We now assume $\forall j\leq i, x_{i+1}\neq x_j$.

Let $h_{i+1}$ the orthogonal projection of $x_{i+1}$ into $\mathcal{H}_i$. 
If $h_{i+1}\neq c_i$, let $v=\frac{(h_{i+1}-c_i)}{\le\|h_{i+1}-c_i\ri\|_2}$. If $h_{i+1}=c_i$, let $s\in\mathcal{S}_{i}$ and $v=\frac{(s-c_i)}{\le\|s-c_i\ri\|_2}$.

Let 
$$c_{i+1}=c_i-\frac{\mu}{L}r_i v$$

$$f_{i+1}=\frac{\mu+L}{4}(\le\|x_{i+1}-c_{i+1}\ri\|_2^2+(1-\frac{\mu^2}{L^2})r_i^2)$$

$$g_{i+1}=L\frac{\le\|x_{i+1}-x^*\ri\|_2}{\le\|x_{i+1}-c_{i+1}\ri\|_2}\le(x_{i+1}-c_{i+1}\ri)$$

$$\mathcal{H}_{i+1}=\le\{x\in\mathcal{H}_i\mid \le\langle x-c_i\mid v\ri\rangle=-\frac{\mu}{L}r_i \ri\}$$

$$\mathcal{S}_{i+1}=\mathcal{S}_i\cap\mathcal{H}_{i+1}$$

$v$ is the difference between two points of $\mathcal{H}_i$, therefore it is one of the direction of $\mathcal{H}_i$, and since $c_i\in\mathcal{H}_i$, $\mathcal{H}_{i+1}$ indeed defines an affine subspace of $\mathcal{H}_i$ of dimension $(d-i-2)$. Let $\mathcal{C}$ the sphere in $\mathcal{H}_{i+1}$ of center $c_{i+1}\in\mathcal{H}_{i+1}$ and radius $r_{i+1}=\sqrt{1-\frac{\mu^2}{L^2}}r_i$. We now want to prove that $\mathcal{C}=\mathcal{S}_{i+1}$.

First, let $x \in \mathcal{H}_{i+1}$. Then
\begin{align}
    \braket{x - c_{i+1}}{v} &= \braket{x - c_i + \frac{\mu}{L} r_i v}{v}    \nonumber \\
   &= \braket{x - c_i}{v} + \frac{\mu}{L} r_i      \nonumber\\
    &= - \frac{\mu}{L} r_i + \frac{\mu}{L} r_i = 0      
    \label{eq_lower1}
\end{align}
\textbf{i)} First, we show that $\mathcal{C} \subseteq \mathcal{S}_{i+1}$.
Let $x\in\mathcal{C}$
\begin{align}
    \norm{x - c_i}^2  &= \norm{x - c_{i+1} -\frac{\mu}{L} r_i v}^2      \nonumber \\
    &= \norm{x - c_{i+1}}^2 + \frac{\mu^2}{L^2} r_i^2 - 2\frac{\mu}{L} r_i\braket{x - c_{i+1}}{v} \nonumber \\
    &= \norm{x - c_{i+1}}^2 + \frac{\mu^2}{L^2} r_i^2       && \text{using \eqref{eq_lower1} since $x\in\mathcal{C}\subseteq\mathcal{H}_{i+1}$} \nonumber \\
    &= (1 - \frac{\mu^2}{L^2}) r_i^2 + \frac{\mu^2}{L^2} r_i^2 \nonumber \\
    &= r_i^2 
\end{align}
since $x\in\mathcal{H}_{i+1}\subseteq\mathcal{H}_i$ and $\le\|x-c_i\ri\|_2=r_i$, $x\in\mathcal{S}_i$ and therefore $x\in\mathcal{S}_{i+1}$

\textbf{ii)} Conversely, we show that $\mathcal{S}_{i+1}\subseteq\mathcal{C}$. Let $x \in \mathcal{S}_{i+1}$. 

\begin{align}
    r_i^2 &= \norm{x - c_i}^2       && \text{as } x \in \mathcal{S}_{i+1}\subseteq\mathcal{S}_i \nonumber \\
    &= \norm{x - c_{i+1} -\frac{\mu}{L} r_i v}^2          \nonumber \\
    &= \norm{x - c_{i+1}}^2 + \frac{\mu^2}{L^2} r_i^2 - 2\frac{\mu}{L} r_i\braket{x - c_{i+1}}{v} \nonumber \\
    &= \norm{x - c_{i+1}}^2 + \frac{\mu^2}{L^2} r_i^2       && \text{using \eqref{eq_lower1} since $x\in\mathcal{S}_{i+1}\subseteq\mathcal{H}_{i+1}$}
\end{align}
from which we obtain that 
\begin{align}
    \norm{x - c_{i+1}}^2 = r_{i+1}^2
\end{align}
So $x\in\mathcal{H}_{i+1}$ and $\le\|x-c_{i+1}\ri\|_2=r_{i+1}$, thus $x\in\mathcal{C}$. We have thus proved that $\mathcal{S}_{i+1}$ is indeed a $(d-i-3)$-sphere in a $(d-i-2)$ affine space with the desired radius and center which concludes the first item of the induction.

\textbf{We now want to prove the second item}. For $x^*\in\mathcal{S}_{i+1}$, 
$x_{i+1}-h_{i+1}$ is orthogonal to $\mathcal{H}_{i+1}$ due to being orthogonal to $\mathcal{H}_i$ by construction. $h_{i+1}-c_{i+1}$ is aligned with $v$ and thus is orthogonal to $\mathcal{H}_{i+1}$. Therefore, their sum $x_{i+1}-c_{i+1}$ is orthogonal to $\mathcal{H}_{i+1}$ and we get
\begin{equation}
\label{proof_low_2}
\begin{aligned}
    \le\|x_{i+1}-x^*\ri\|_2^2 & = \le\|x_{i+1}-c_{i+1}\ri\|_2^2+\le\|c_{i+1}-x^*\ri\|_2^2 \\
    & =\le\|x_{i+1}-c_{i+1}\ri\|_2^2+r_{i+1}^2 \\
\end{aligned}
\end{equation}
And thus 
\begin{equation}
\label{cond1}
    f_{i+1}=\frac{\mu+L}{4}\le\|x_{i+1}-x^*\ri\|_2^2
\end{equation}

Let $x^*\in\mathcal{S}_{i+1}$. Since $\mathcal{S}_{i+1}\subseteq \mathcal{S}_{i}$, then by recurrence hypothesis there exists an interpolation of the $(x_j,f_j,g_j)$ in $RSI^-(\mu)\cap EB^+(L)$ minimized by $x^*$, hence from Theorem \ref{interpolation_th}, 
\begin{equation} 
\label{cond3}
\forall j\leq i, \le\|g_{j}\ri\|_2\leq L\le\|x_j-x^*\ri\|_2 \quad \text{and}\quad \le\langle g_{j}\mid x_j-x^*\ri\rangle\geq\mu\le\|x_j-x^*\ri\|_2^2
\end{equation}
Moreover, by construction of $g_{i+1}$,
\begin{equation}
\label{cond2}
    \le\|g_{i+1}\ri\|_2=L\le\|x_{i+1}-x^*\ri\|_2
\end{equation}

Since $x_{i+1}-c_{i+1}$ is orthogonal to $\mathcal{H}_{i+1}$ and thus to $c_{i+1}-x^*$, we have 
\begin{equation} 
\label{prereq1}
\le\langle x_{i+1}-c_{i+1}\mid x_{i+1}-x^*\ri\rangle=\le\|c_{i+1}-x^*\ri\|_2^2
\end{equation} 
Besides, $x_{i+1}-c_{i+1}$ is orthogonal to $c_{i+1}-x^*$ and thus
\begin{equation}
    \label{prereq2}
\le\|x_{i+1}-x^*\ri\|_2^2=\le\|x_{i+1}-c_{i+1}+c_{i+1}-x^*\ri\|_2^2=\le\|x_{i+1}-c_{i+1}\ri\|_2^2+r_{i+1}^2
\end{equation}

By construction,

\begin{equation}
    \label{prereq3}
    \begin{aligned}
    \le\|c_{i+1}-x_{i+1}\ri\|_2 & \geq \le\|c_{i+1}-h_{i+1}\ri\|_2 \\
    & = \le\|c_{i+1}-c_i+c_i-h_{i+1}\ri\|_2 \\
    & = \le\|-\frac{\mu}{L}r_i v - \le\|h_{i+1}-c_i\ri\|_2 v\ri\|_2 \\
    & = \frac{\mu}{L}r_i + \le\|h_{i+1}-c_i\ri\|_2 \\
    & \geq \frac{\mu}{L}r_i \\
    \end{aligned}
\end{equation}

and finally :

\begin{equation}
\begin{aligned}
\frac{\le\langle g_{i+1}\mid x_{i+1}-x^*\ri\rangle^2}{\mu^2\le\|x_{i+1}-x^*\ri\|_2^4} & =  \frac{L^2}{\mu^2}\frac{\le\langle x_{i+1}-c_{i+1}\mid x_{i+1}-x^*\ri\rangle^2}{\le\|x_{i+1}-x^*\ri\|_2^2\le\|x_{i+1}-c_{i+1}\ri\|_2^2} \\
& = \frac{L^2}{\mu^2}\frac{\le\|x_{i+1}-c_{i+1}\ri\|_2^2}{\le\|x_{i+1}-x^*\ri\|_2^2} && \text{using \eqref{prereq1}}\\
& = \frac{L^2}{\mu^2}\frac{\le\|x_{i+1}-c_{i+1}\ri\|_2^2}{\le\|x_{i+1}-c_{i+1}\ri\|_2^2+r_{i+1}^2} && \text{using \eqref{prereq2}}\\
& \geq \frac{L^2}{\mu^2}\frac{\frac{\mu^2}{L^2}r_i^2}{\frac{\mu^2}{L^2}r_i^2+(1-\frac{\mu^2}{L^2})r_i^2} && \text{using \eqref{prereq3}}\\
& = 1
\end{aligned}
\end{equation}

And thus 

\begin{equation}
    \label{cond4}
    \le\langle g_{i+1}\mid x_{i+1}-x^*\ri\rangle\geq\mu\le\|x_{i+1}-x^*\ri\|_2^2
\end{equation}

Since $\forall j\leq i, x^*\in\mathcal{S}_j$, we also have 
\begin{equation}
    \label{cond5}
    f_j=\frac{\mu+L}{4}\le\|x_j-x^*\ri\|_2^2
\end{equation}
We finally apply Corollary \ref{interp_cor} to all \textit{unique} triples $(x_j, f_j, g_j)_{j\leq i+1}$ (which ensures by construction that all $x_j$ are distincts) which allows us to conclude from \eqref{cond1}, \eqref{cond2}, \eqref{cond3}, \eqref{cond4} and \eqref{cond5} that there exists an interpolation in $RSI^-(\mu)\cap EB^+(L)$ that is minimized by $\{x^*\}$, proving the second item of the induction and thus concluding the proof.

$\hfill\blacksquare$

\end{document}